\documentclass[11pt]{article}

\usepackage{amsfonts,amsthm,amsmath,amsfonts,latexsym,amssymb}

\usepackage{graphicx,upref,calligra}
\usepackage[svgnames]{xcolor}
\usepackage{mathrsfs}
\usepackage{authblk}
\usepackage[T1]{fontenc}
\usepackage{enumerate}
\usepackage[all]{xy}
\usepackage{bbm}

\newtheorem{teo}{Theorem}[section]

\newtheorem{pro}[teo]{Proposition}
\newtheorem{lem}[teo]{Lemma}
\newtheorem{cor}[teo]{Corollary}
\newtheorem{defi}[teo]{Definition}
\theoremstyle{definition}
\newtheorem{rem}[teo]{Remark}
\newtheorem{exa}[teo]{Example}

\newtheorem{nota}[teo]{Notation}
\def\QED{\hfill $\blacksquare$}
\def\EOE{\hfill $\blacktriangle$}

\def\bdem{\begin{proof}}
\def\edem{\end{proof}}
\newcommand{\peso}[1]{ \ \ \text{ \rm  #1 } \ \ }
\def\noi{\noindent}

\definecolor{myblue}{rgb}{0,0.33,0.55}
\definecolor{azul}{rgb}{0.1,0.6,0.86}

\definecolor{myyellow}{rgb}{0.42,0.24,0.52}
\definecolor{mygreen}{rgb}{0.12,0.5,0.29}
\definecolor{myred}{rgb}{0.74,0.13,0.13}
\definecolor{mylgreen}{rgb}{0.68,0.98,0.6}
\definecolor{mylyellow}{rgb}{0.86,0.85,0.55}
\definecolor{myllyellow}{rgb}{0.87,0.86,0.56}
\definecolor{naranja}{RGB}{249,153,96}

\def\la{\lambda}

\def\uno{\mathbbm{1}} 

\def\N{\mathbb{N}}   
\def\R{\mathbb{R}}   
\def\C{\mathbb{C}}   

\def\T{\mathbb{T}}
\DeclareMathOperator{\Span}{Span}

\def\cC{\mathcal{C}}

\def\cG{\mathcal{G}}

\def\cH{\mathcal{H}}

\def\cP{\mathcal{P}}

\def\cM{{\cal M}}

\def\ese{\mathcal{S}}

\def\eme{\mathcal{M}}

\def\cU{\mathcal{U}}

\def\beq{\begin{equation}}
\def\eeq{\end{equation}}

\def\rk{\text{\rm rk}}
\def\orto{^\perp}
\def\coma{\, , \, }
\def\pausa{\medskip\noi}

\def\da{^\downarrow}
\def\bm{\left[\begin{array}}
\def\em{\end{array}\right]}
\def\ben{\begin{enumerate}}
\def\een{\end{enumerate}}
\def\bit{\begin{itemize}}
\def\eit{\end{itemize}}
\def\barr{\begin{array}}
\def\earr{\end{array}}
\def\igdef{\ \stackrel{\mbox{\tiny{def}}}{=}\ }

\newcommand{\diag}[1]{\hbox{diag}\left( #1\right)}

 \DeclareMathOperator{\tr}{tr}

\newcommand{\pint}[1]{\displaystyle \left \langle\, #1 \, \right\rangle}

\newcommand{\hil}{\mathcal{H}}

\newcommand{\mat}{\mathcal{M}_n(\mathbb{C})}

\newcommand{\matsa}{\mathcal{H}(n)}
\newcommand{\matsai}{\mathcal{H}_n(I)}

\newcommand{\matu}{\mathcal{U}(n)}
\newcommand{\matpos}{\mat^+}
\newcommand{\matppos}{\matinv^+}

\newcommand{\matinv}{\mathcal{G}\textit{l}\,(n)}

\newcommand{\spec}[1]{\sigma\left(#1\right)}

\newcommand{\avi}[2]{\la_{#1}\left( #2\right)}
\newcommand{\svi}[2]{s_{#1}\left( #2\right)}

\newcommand{\sub}[2]{{#1}_{\mbox{\tiny{${#2}$}}}}

\newcommand{\convm}{\xrightarrow[m\rightarrow\infty]{}}

\newcommand{\mnorm}[1]{%
\left\vert\kern-0.9pt\left\vert\kern-0.9pt\left\vert #1
\right\vert\kern-0.9pt\right\vert\kern-0.9pt\right\vert}
\newcommand{\bmnorm}[1]{%
\big\vert\kern-0.9pt\big\vert\kern-0.9pt\big\vert #1
\big\vert\kern-0.9pt\big\vert\kern-0.9pt\big\vert}

\begin{document}

\title{{\bf{Minimal curves in $\matu$ and $\matppos$ with respect to the spectral and the trace norms}}}

\author[1,2]{Jorge Antezana}
\author[1, 2]{Eduardo Ghiglioni}
\author[1, 2]{Demetrio Stojanoff}
\affil[1]{Departamento de de Matem\'atica, FCE-UNLP, Calles 50 y 115, 
	(1900) La Plata, Argentina.}
\affil[2]{ Instituto Argentino de Matem\'atica, `Alberto P. Calder\'on', CONICET, Saavedra 15 3er. piso,
	(1083) Buenos Aires, Argentina.}

\date{}

\maketitle

\begin{abstract}
Consider the Lie group of $n\times n$ complex unitary matrices $\matu$ endowed with the bi-invariant Finsler metric given by the spectral norm,
$$
\sub{\|X\|}{U}=\|U^*X\|_{\infty}=\|X\|_{\infty}
$$
for any $X$ tangent to a unitary operator $U$. Given two points in $\matu$, in general there exists infinitely many curves of minimal length. The aim of this paper is to provide a complete description of such curves. As a consequence of this description, we conclude that there is a unique curve of minimal length between $U$ and $V$ if and only if the spectrum of $U^*V$ is contained in a set of the form $\{e^{i\theta},e^{-i\theta}\}$ for some $\theta\in [0,\pi)$. Similar
studies are done for the Grassmann manifolds.
\\
\indent
Now consider the cone of $n \times n$ positive invertible matrices $\matppos$
endowed with the bi-invariant Finsler metric given by the trace norm,
$$
\sub{\|X\|}{1,A}=\|A^{-1/2}XA^{-1/2}\|_1
$$
for any $X$ tangent to $A\in\matppos$. In this context, given two points $A, B \in \matppos$ there exists infinitely many curves of minimal length. In order to provide a complete description of such curves, we provide a characterization of the minimal curves joining two Hermitian matrices $X,Y\in \matsa$. As a consequence of the last description, we provide a way to construct minimal paths in the group of unitary matrices $\matu$ endowed with the bi-invariant Finsler metric
$$
\sub{\|X\|}{1,U}=\|U^*X\|_1=\|X\|_1
$$
for any $X$ tangent to $U\in\matu$.
\\
\indent
 We also study the set of intermediate points in all the previous contexts. Between two given unitary matrices $U$ and $V$ we prove that this set is geodesically convex provided $\|U-V\|_{\infty}<1$. In $\matppos$ this set is geodesically convex for every unitarily invariant norm.
\end{abstract}

\tableofcontents

\section{Introduction}

Let $\matinv$ denote the Lie group of invertible $n\times n$ matrices, $\matu$ the Lie subgroup of unitary matrices, and $ \matppos$ the cone of positive invertible matrices. Given  $T\in\matinv$, it can be decomposed as
$$
T=U|T|,
$$
where $U$ is a unitary matrix, and $|T|$ is the positive matrix given by $(T^*T)^{1/2}$. This is the usual polar decomposition, and a particular case of the so called Cartan decomposition for more general Lie groups. This decomposition allows to understand the metric and geometric properties of $\matinv$ through the  study of geometric and metric properties of $\matu$ and $\matppos$.

\medskip

In the case of $\matu$, as in any Lie group, it has a canonical torsion-free connection defined on left-invariant vector fields $X,Y$ by $\nabla_XY=\frac12 [X,Y]$. The geodesics associated to this connection are the one-parameter groups  $t\mapsto Ue^{tZ}$ (here $U$ is a unitary matrix and $Z$ an anti-hermitian matrix). We can also introduce a Riemannian metric on the unitary group in a standard way
$$
\langle X,Y\rangle_U=Tr(U^*X(U^*Y)^*)=Tr(XY^*),
$$
for $U^*X,U^*Y$ in the Lie algebra of the group, that is, for $U^*X,U^*Y$ anti-Hermitian matrices. It is well-known that the aforementioned connection is the Levi-Civita connection of the metric induced by the trace, and that geodesics are minimal curves for $t\in[0,1]$ provided the spectrum of $iZ$ is contained in $(-\pi,\pi)$. 
 
\medskip

With respect to $\matppos$,  it is an open subset of the space of hermitian matrices $\matsa$. Therefore, it inheres a geometric structure where the tangent spaces can be identified with $\matsa$. Also, there exists a natural transitive action of $\matinv$ by conjugation. The properties of this action make $\matppos$ become an homogeneous space. Moreover, using this action it is possible to define a covariant derivative that leads to the following differential equation for the geodesics
$$
\gamma''=\gamma'\gamma^{-1}\gamma',
$$
 (see \cite{CoMa}, \cite{CPR3} and \cite{KN}). Given $A,B\in\matppos$, the solution of the corresponding Dirichlet problem gives the following expression for the geodesic joining $A$ with $B$
 $$
\sub{\gamma}{\matppos}(t)=A^{1/2}(A^{-1/2}BA^{-1/2})^tA^{1/2}.
 $$  
Note that if one of the end points is the identity, also in this case the geodesic is a one-parameter subgroup of $\matinv$. 

\medskip
It is also possible to define  a (canonical) Riemannian structure on $\matppos$. Indeed, consider the inner product associated to the trace in the tangent space of $\matppos$ at the identity. Then, using the homogeneous structure, we can define the following inner product
$$
\pint{X,Y}_A=\tr\big((A^{-1/2}XA^{-1/2})(A^{-1/2}YA^{-1/2})^*\big)=\tr(A^{-1/2}XA^{-1}Y^*A^{-1/2}),
$$
in the tangent space corresponding to another  point $A\in\matppos$. Endowed with this structure, the action by conjugations becomes isometric, and $\matppos$ becomes a complete Riemannian manifold with non-positive sectional curvature. As before, these geodesics are minimal curves. However, in this case, they are short not only for $t\in [0,1]$, but also for every $t\in\R$  (see \cite{Bh PM}, \cite{PR}).

\bigskip

Now, let $\|\cdot\|_p$ denote the Schatten norm defined by
\begin{align*}
\|A\|_p&=
\tr(|A|^p)^{1/p}\quad \ \ \ \  \ \ \ \ \mbox{if $1\leq p<\infty$},\\
\|A\|_{\infty}&=\sup_{x\in\C^n\setminus\{0\}}\frac{\|Ax\|}{\|x\|}.
\end{align*}
Using these norms, we can define the bi-invariant Finsler metric such that
$$
\sub{\|X\|}{p,U}=\|U^*X\|_p=\|X\|_p
$$
in the tangent space at $U\in\matu$ and
$$
\sub{\|X\|}{p,A}=\|A^{-1/2}XA^{-1/2}\|_p
$$
in the tangent space at $A\in\matppos$. With respect to any of these Finsler metrics, the length of a curve $\alpha$ parametrized by the interval $[a,b]$ is computed by 
$$
L(\alpha)=\int_a^b \|\dot{\alpha}\|_{p,\alpha}\,dt.
$$
The induced rectifiable distance, denoted by $d_p(\cdot,\cdot)$, is computed as in the Riemannian setting as the infimum of the length of piecewise smooth curves joining given endpoints  (see \cite{ALV} and the references therein). Some applications of these metrics to control theory have been found in \cite{ALV}.

\medskip

A remarkable fact is that the aforementioned geodesics remain minimal curves with respect to all these new metric structures. Moreover, if $1<p<\infty$, they are the unique minimal curves joining two elements (see  \cite{An05} \cite{upe}, \cite{AR06tesisChumiento}, \cite{AR08tesisChumiento}, \cite{ALV} for the unitary case and  \cite{bhatiapositivas},  \cite{DMR04}, \cite{DMR05}, \cite{langpositivas}, \cite{lawsonpositivas} for the positive case). As in the Riemannian case, in the case of the group of unitary matrices, the uniqueness holds for $t\in [0,1]$ only if $\|Z\|_\infty<\pi$. 

\medskip

If $p=\infty$ or $p=1$, the situation is different, and there exist infinitely many minimal curves joining two points either in $\matu$ or in $\matppos$. In \cite{L1} Yongdo Lim study the set of minimal curves in $\matppos$ endowed with the Finsler structure associated to $\|\cdot\|_\infty$. He completely characterized all the minimal curves. Moreover, he also studied the sets of midpoints
 $$
\eme_{1/2}(A,B)=\big\{C\in\matppos: d_\infty(A,C)=d_\infty(C,B)= \mbox{$\frac{1}{2}$}\,d_\infty(A,B) \big\},
$$
or more generally for $t\in(0,1)$ the set of intermediate points
$$
\eme_t(A,B)=\big\{C\in\matppos: d_\infty(A,C)=t\,d_\infty(A,B),\  d_\infty(C,B)= (1-t)\,d_\infty(A,B)\big\}.
$$

The main aim of this work is to continue Lim's study in the remaining cases, that is, to characterize the minimal curves in the following cases:
\begin{itemize}
\item The space $\matppos$ with the Finsler structure associated to $\|\cdot\|_1$;
\item The Lie group $\matu$ with the Finsler structures associated to $\|\cdot\|_\infty$ and $\|\cdot\|_1$.
\end{itemize}

 The strategies to get a characterization of the minimal curves is different in each case. In the case of unitary matrices with $p=\infty$, the key results are Lemma \ref{CPR1} and Lemma \ref{CPR2} which are based on a beautiful trick used by Porta and Recht in \cite{PR}. As a consequence of our characterization, we get a description of those pairs for which there exists a unique minimal curve. In this case there is a unique minimal curve connecting $U,V\in\matu$ if and only if the spectrum of $U^*V$ is contained in $\{e^{i\theta},e^{-i\theta}\}$ for some $\theta\in [0,\pi)$. This is similar to the result obtained by Lim in \cite{L1} in the case of positive operators. On the other hand, the Grassmannian can be modeled as a submanifold of the unitary group
(identifying a subspace with the associated orthogonal symmetry). Using this 
idea we also describe all minimal curves connecting two projections $P$ and $Q$ 
such that  $\left\|P - Q\right\|_{\infty} < 1$. 

\medskip

In the case of $\matppos$ and $p=1$, we get the characterization by lifting the problem to the space of hermitian matrices $\matsa$. So, firstly we provide a characterization of the minimal paths joining $X,Y\in \matsa$, if the length of a curve $\alpha:[a,b]\to\matsa$ is measure by
$$
L(\alpha)=\int_a^b \|\dot{\alpha}\|_1\,dt.
$$
For the best of our knowledge, this characterization is not known. Our argument is based in a trick that reduces the problem to the $2\times 2$ case (see Lemma \ref{dospordos} and Theorem \ref{caracterizacion}). Once the characterization of the minimal curves is given for $\matsa$, the characterization in $\matppos$ can be obtained using the so called Exponential Metric Increasing property (see \cite{Bh PM}, and also \cite{CPR2}).

\bigskip

Finally, in the case of $\matu$ and $p=1$ the above lifting argument also works in one direction. Indeed, using the same idea as in the case of $\matppos$, we prove that a minimal curve in $\matsa$ leads to a minimal curve in $\matu$ by means of the exponential map. However, the corresponding EMI property can not be used in this case to prove the converse. Roughly speaking, the reason is that in this case the space has positive curvature, and the corresponding exponential inequality reverses its direction. However, our characterization in $\matsa$ allow us to construct minimal curves in $\matu$.

\medskip

Also, to continue Lim's study \cite{L1} we prove that the set of intermediate points is geodesically convex in all the previous contexts.  Actually, we will use the same idea for all this cases. First we prove that the function
$$
s \mapsto d(I, \gamma(s)),
$$
where $\gamma(s)$ is the geodesic joining two matrices (in some of this spaces for some metric), is convex. Then we will use this fact to prove that the set of intermediate points is geodesically convex. Indeed, between two given unitary matrices $U$ and $V$ we prove that it is geodesically convex provided $\|U-V\|_{\infty}<1$. In $\matppos$ it's geodesically convex for every unitarily invariant norm. We note that the set of intermediate points in $\matppos$ was already proved to be geodesically convex for every unitarily invariant norm in \cite{L1} but using a different technique.

\bigskip

The paper is organized as follows:  in section 2 we describe all the possible short paths connecting two points in $\matu$. In section 3 we describe all the possible short paths connecting two points in the Grassmannian, and in section 4 the case of the trace norm is studied. First we get the characterization of minimal curves in the space of Hermitian matrices $\matsa$ and after that we study the space $\matppos$. Finally, in section 5 we study the geometry of intermediate points.

\subsection*{Notation}
 
Throughout this note, $\mat$ denotes the algebra of complex $n\times n$ matrices, $\matinv \subseteq \mat$ the group of all invertible matrices,  $\matu$ the group of unitary $n\times n$ matrices, and $\matsa$ the real subspace of Hermitian matrices.  Sometimes we will write $\cG l(\ese)$ (resp. $\cU(\ese)$, $\hil(\ese)$) to indicate that the operators are acting on some specific subspace $\ese$. 

\pausa
If $T\in \mat$, then $\|T\|_{\infty}$ stands for the usual spectral norm,  $|T|$  indicates the modulus of $T$, i.e. $|T|=\sqrt{T^*T}$, and $\tr(T)$ denotes the trace of $T$. By means of $\spec{T}$ we denote the set of eigenvalues of $T$, while $\rho(T)$ denotes the spectral radius of $T$. 

\pausa
Given $A\in\matsa$, $\avi{1}{A}\geq \ldots\geq \avi{n}{A}$ denotes the eigenvalues of $A$ arranged in non-increasing way and counted with multiplicity. Analogously,  given an arbitrary matrix $T\in\mat$, $\svi{1}{T}\geq \ldots\geq \svi{n}{T}$ denotes the singular values of $T$ (also counted with multiplicity), i.e. the eigenvalues of $|T|$. Finally, given $A,B\in\matsa$, by means of $A\leq B$ we denote that $A$ is less that or equal to $B$ with respect to the L\"owner order.

\section{The spectral norm}

In this section, we will study the structure of minimal curves joining two given unitary matrices $U$ and $V$, where the minimality is with respect to the Finsler metric structure given by the spectral norm. 

\paragraph{Convention:} From now on, we assume that the curves are parametrized by the interval $[0,1]$, and  in such a way that $t\mapsto\|\dot\alpha(t)\|$ is constant.

\bigskip
\subsection{Structure of minimal curves}
Recall that, in the case of hermitian matrices,  the spectral norm coincides with the spectral radius. Thus, roughly speaking, we can change an hermitian matrix inside the eigenspace corresponding to the ``small eigenvalues'' and it will still have the same spectral norm. This simple observation gives an easy strategy to construct many minimal curves by perturbing the uniparametric groups. 

\medskip

The following result says that precisely these are all the possible minimal curves. For the sake of simplicity, and without lost of generality, we will assume that one of the endpoints is the identity.

\medskip

\begin{teo}\label{structure}
Given $U\in\matu$, let $X\in\matsa$ such that $U=e^{iX}$ and $\|X\|_{\infty}\leq \pi$. Then:
\begin{enumerate}
\item[a)] If $\spec{|X|}\subseteq [0,\pi)$ and it has more than one element, then the minimal curves joining $I$ with $U$ have the following structure:
$$
\alpha(t)=\begin{pmatrix}
e^{it\sub{X}{\ese}}&0\\
0&\alpha_1(t)
\end{pmatrix}\begin{array}{l}
\ese\\
\ese^\bot
\end{array},
$$
where $\ese=\ker(\|X\|_{\infty}I-|X|)$, $\sub{X}{\ese}=\left.X\right|_{\ese}\in \cH(\ese)$, and $\alpha_1:[0,1]\to \cU(\ese^\bot)$ is any curve joining $\sub{I}{\ese^\bot}$ and $\sub{U}{\ese^\bot}=\left.U\right|_{\ese^\bot}$ such that $\|\dot\alpha_1\|_{\infty}\leq \|X\|_{\infty}$.
\item[b)] If $\pi\in \spec{|X|}$, then the minimal curves joining $I$ with $U$ have the following structure:
$$
\alpha(t)=\begin{pmatrix}
e^{itY_\ese}&0\\
0&\alpha_1(t)
\end{pmatrix}\begin{array}{l}
\ese\\
\ese^\bot
\end{array},
$$
where $\ese=\ker(\pi I-|X|)$, $Y_\ese\in \hil(\ese)$ and it satisfies that $\spec{Y_\ese}\subseteq \{\pi, -\pi\}$, and $\alpha_1:[0,1]\to \cU(\ese^\bot)$ is any curve joining $\sub{I}{\ese^\bot}$ and $\sub{U}{\ese^\bot}=\left.U\right|_{\ese^\bot}$ such that $\|\dot\alpha_1\|_{\infty}\leq \pi$.
\end{enumerate}
\end{teo}

A direct consequence of the first item of this result is the following corollary about the uniqueness of minimal curves:

\begin{cor}\label{unica}
Given $U,V\in\matu$, the exists a unique minimal curve between them if and only if $\spec{U^*V}\subseteq\{e^{ i\theta}, e^{ -i\theta}\}$ for some $|\theta|<\pi$.
\end{cor}

Another consequence of this theorem is the following result:

\begin{cor}\label{lineal}
Given a minimal curve $\alpha:[0,1]\to\matu$, it is also minimal between any two points of its trace, and 
$$
d_\infty(U,\alpha(r))=rd_\infty(U,V).
$$
\end{cor}

\subsection{Proof of Theorem \ref{structure}}

To begin with, we will prove the following lemma, which is a modification of a beautiful trick used by Porta and Recht in \cite{PR} (see also \cite{ALV}):

\begin{lem}\label{CPR1}
Let $\alpha:[0,1]\to\matu$ be a (smooth) minimal curve joining $I$ with $U=e^{iX}$, where $\|X\|_{\infty}<\pi$. If $\xi$ is a unitary eigenvector of $X$ corresponding to an eigenvalue $\la$ such that $\rho(X)=|\la|$. Then
$$
\alpha(t)=\begin{pmatrix}
e^{it\la}&0\\
0&\alpha_1(t)
\end{pmatrix}\begin{array}{l}
\ese\\
\ese^\bot
\end{array},
$$
where $\ese=\mbox{span}\{\xi\}$, and $\alpha_1:[0,1]\to \cU(\ese^\bot)$ is any curve joining $\sub{I}{\ese^\bot}$ and $\sub{U}{\ese^\bot}=\left.U\right|_{\ese^\bot}$ such that $\|\dot\alpha_1\|_{\infty}\leq \|X\|_{\infty}$.
\end{lem}

\medskip

\bdem
Let $\gamma:[0,1]\to\matu$ be the curve defined by $\gamma(t)=e^{itX}$, and let $S^{2n-1}$ denote the sphere:
$$
S^{2n-1}=\{\eta\in\C^n: \|\eta\|_2=1\}.
$$
Define the curves $\widetilde{\gamma}:[0,1]\to S^{2n-1}$ and  $\widetilde{\alpha}:[0,1]\to S^{2n-1}$ in the following way:
$$
\widetilde{\gamma}(t)=\gamma(t)\xi\peso{and}\widetilde{\alpha}(t)=\alpha(t)\xi.
$$

A simple computation shows that 
$
\mbox{L}(\widetilde{\gamma})=\mbox{L}(\widetilde{\alpha})
$, 
where the length here is computed with respect to the Riemannian structure of the sphere. Since $\|X\|_{\infty}<\pi$, $\widetilde{\gamma}$ is the unique geodesic in $S^{2n-1}$ joining $\xi$ with $U\xi=e^{iX}\xi$. Hence, $\widetilde{\gamma}=\widetilde{\alpha}$. In particular, $\xi$ is an eigenvector of $\alpha(t)$ associated to $e^{it\la}$ for every $t\in[0,1]$. The rest of the statement is a consequence of the block decomposition of $\alpha$, induced by the decomposition of the space $\C^n=\ese\oplus\ese^\bot$.
\edem

Iterating this lemma we get the first part of Theorem \ref{structure}.

\begin{cor}\label{primera parte}
Let $\alpha:[0,1]\to\matu$ be a (smooth) minimal curve joining $I$ with $U=e^{iX}$, where $\|X\|_{\infty}<\pi$. Then 
$$
\alpha(t)=\begin{pmatrix}
e^{it\sub{X}{\ese}}&0\\
0&\alpha_1(t)
\end{pmatrix}\begin{array}{l}
\ese\\
\ese^\bot
\end{array},
$$
where $\ese=\ker(\|X\|_{\infty}I-|X|)$, $\sub{X}{\ese}=\left.X\right|_{\ese}$, and $\alpha_1:[0,1]\to \cU(\ese^\bot)$ is a curve joining $\sub{I}{\ese^\bot}$ and $\sub{U}{\ese^\bot}=\left.U\right|_{\ese^\bot}$ such that $\|\dot\alpha_1\|_{\infty}\leq \|X\|_{\infty}$.
\end{cor}

For the second part, we will use an approximation argument.

\begin{lem}\label{CPR2}
Let $\alpha:[0,1]\to\matu$ be a (smooth) minimal curve joining $I$ with $U=e^{iX}$, where $\|X\|_{\infty}=\pi$. There exists $Y \in \mathcal{H}(n)$ with $\left\|Y\right\|_{\infty} = \pi$ such that $U = e^{iY}$ and a unitary eigenvector $\xi$ of $Y$ corresponding to an eigenvalue $\la$ with $|\la|=\pi$ such that
$$
\alpha(t)=\begin{pmatrix}
e^{it\la}&0\\
0&\alpha_1(t)
\end{pmatrix}\begin{array}{l}
\ese\\
\ese^\bot
\end{array},
$$
where $\ese=\mbox{span}\{\xi\}$, and $\alpha_1:[0,1]\to \cU(\ese^\bot)$ is any curve joining $\sub{I}{\ese^\bot}$ and $\sub{U}{\ese^\bot}=\left.U\right|_{\ese^\bot}$ such that $\|\dot\alpha_1\|_{\infty}\leq \|X\|_{\infty}$.
\end{lem}

\bdem
For each $ m\in \N$ we will consider the curve 
$$
\alpha_{m}(t) = \alpha\left(\left(\frac{m}{m + 1}\right)t\right) \ ,
\peso{for every} t \in [0, 1] \ .
$$
In other words, $\alpha_{m}(t)$ is the part of the curve $\alpha(t)$, joining $I$ with $\alpha_{m}(1) = e^{iX_{m}}$ for some $X_{m} \in \matsa$ with $\|X_m\|_{\infty} < \pi$. As $\alpha(t)$ has constant speed, it holds that
$$
\left\|X_{m}\right\|_{\infty} = d_\infty(I,e^{iX_m})=\frac{m}{m + 1} L(\alpha) < \pi. 
$$
Therefore, $\alpha_m(t)$ is a minimal curve joining $I$ with $e^{iX_m}$. Moreover, 
$$
e^{iX_{m}} = \alpha \left(\frac{m}{m + 1}\right) \xrightarrow[m \rightarrow \infty]{}  \alpha(1) = e^{iX}.
$$
For $\left\{X_m\right\}_{m \in \N}$ there is a convergent subsequence, which we denote $X_m$ with some abuse of notation, i.e, 
\begin{equation}\label{convergenciaoperadores}
X_m \xrightarrow[m \rightarrow \infty]{} Y \in \matsa \implies U = e^{iX}= e^{iY}\ .
\end{equation}
Denote by $U_m = e^{iX_{m}}$ for $m\in\N$, and 
$$
\eta := \min \left\{\left|\lambda_j(U) - \lambda_k(U)\right| : \ \mbox{where} \ \lambda_j(U) \neq \lambda_k(U)\right\} > 0 \ .
$$
For each $j$ consider the disc $D_j$ with center $\lambda_j(U)$ and radius $\eta/3$. By \eqref{convergenciaoperadores} and the fact that 
$\sigma(U) \subseteq \bigcup_{j} D_j\,$, there is an $N_1 \geq 1$ such that
$$
\sigma(U_m) \subseteq \bigcup_{j} D_j \peso{for every} m \geq N_1\ .
$$
Let $D_0 = \overline{D}(-1, \eta/3)$ and 
take a continuous map $f : \mathbb{T} \rightarrow [0, 1]$ such that 
$f |_{D_0} = 1$ and $ f |_{D_k} = 0$ for $k \neq 0$. 
Then $f(U_m)\ \xrightarrow[m \rightarrow \infty]{} f(U)$ by the Spectral Theorem. 
In other words
\begin{equation}\label{proyectores}
P_{\mathcal{T}_m} \xrightarrow[m \rightarrow \infty]{} P_{\mathcal{T}}\ ,
\end{equation}
where $P_\mathcal{T}= f(U)$ is the projection onto the subspace 
$\mathcal{T} = \ker(U + I)$ and $P_{\mathcal{T}_m} =f(U_m)$ is the 
projection onto the subspace 
$$
\mathcal{T}_m = \Span \left\{v\in \C^n  : U_m v = \lambda_{k}(U_m)v \peso{for some} \lambda_{k}(U_m) \in D_0\right\}.$$
For each $m\in \N$, let $\la_m  \in\sigma(X_m)$ 
such that $|\la_m|=\left\|X_m\right\|_{\infty}$. As $\la_m \in [-\pi, \pi]$ 
there is a convergent subsequence, which we denote $\la_m$ with 
some abuse of notation again, i.e, 
$$ \la_m  \xrightarrow[m \rightarrow \infty]{} \la, 
$$ 
where $|\la|= \pi$. 
As $\left\|X_m\right\|_{\infty} < \pi$ we can apply the Lemma \ref{CPR1} for the curve $\alpha_m(t)$. 
For each $m\in \N$, let $\xi_m \in \ker (X_m - \la_m \,I)$ be a unit vector, and 
denote by $\ese_m=\mbox{span}\{\xi_m\}$. Then 
$$
\alpha_m(t)=\begin{pmatrix}
e^{it\la_m}&0\\
0&\alpha_{m,1}(t)
\end{pmatrix}\begin{array}{l}
\ese_m\\
\ese_m^\bot
\end{array},
$$
where 
$\alpha_{m,1}:[0,1]\to \cU(\ese_m^\bot)$ is a 
curve joining $\sub{I}{\ese_m^\bot}$ and $\sub{U}{m, \ese_m^\bot}$ such that $\|\dot\alpha_{m,1}\|_{\infty}\leq \|X_m\|_{\infty}\,$. As $\xi_m \in \ese_m \subseteq \mathcal{T}_m$ and 
$P_{\mathcal{T}_m} \xrightarrow[m \rightarrow \infty]{} P_{\mathcal{T}}\,$, 
we conclude that there exists a unitary eigenvector $\xi \in \mathcal{T}$ and a convergent subsequence such that
$$
\xi_{j_m} \xrightarrow[m \rightarrow \infty]{} \xi
\implies 
P_{\mathcal{S}_{j_m}}  \xrightarrow[m \rightarrow \infty]{} P_{\mathcal{S}}  \ ,
$$ 
where $\ese=\mbox{span}\{\xi\}$. 
Since $X_{j_m} \convm Y$ and $\la_{j_m} \convm \la$, then $Y\,\xi = \la\, \xi$. Therefore
$$
\left\langle e^{it\la_{j_m}}\xi_{j_m}, \xi_{j_m} \right\rangle  
\xrightarrow[m \rightarrow \infty]{}  \left\langle e^{it\la}\xi,\xi\right\rangle 
\peso{for every} t\in [0\coma 1] \ .
$$
Since also $\alpha_{j_m}(t)\xrightarrow[m \rightarrow \infty]{} \alpha (t) $ 
for every $t\in [0\coma 1]$, it follows that 
$$
\alpha_{j_m}(t)=\begin{pmatrix}
e^{it\la_{j_m}}&0\\
0&\alpha_{{j_m},1}(t)
\end{pmatrix}\begin{array}{l}
\ese_{j_m}\\
\ese_{j_m}^\bot
\end{array} \xrightarrow[m \rightarrow \infty]{}
\alpha(t)=\begin{pmatrix}
e^{it\la}&0\\
0&\alpha_{1}(t)
\end{pmatrix}\begin{array}{l}
\ese\\
\ese^\bot
\end{array},
$$
where $\alpha_1 = P_{\ese\orto } \, \alpha \, P_{\ese\orto } \big|_{\ese\orto }  :[0,1]\to \cU(\ese^\bot)$ is a curve joining $\sub{I}{\ese^\bot}$ and $\sub{U}{\ese^\bot}=\left.U\right|_{\ese^\bot}$ such that 
$\|\dot\alpha_1\|_{\infty} \le\|\dot\alpha\|_{\infty} = \|X\|_{\infty}$. 
\edem

\bigskip

To complete the proof of Theorem \ref{structure} part b) we will iterate this lemma and, in every step, we will get a new $Y$ which could be different from $X$. The reason of this is that there isn't uniqueness in the geodesics joining two antipodal points in the sphere.    

\bdem[Proof of Theorem \ref{structure}]
The first part has been already proved in Corollary \ref{primera parte}. The second part is a consequence of Lemma \ref{CPR2}:
Let $\alpha:[0,1]\to\matu$ be a (smooth) minimal curve joining $I$ with $U=e^{iX}$, where $\|X\|_{\infty}=\pi$. By Lemma \ref{CPR2} there exists $Y_1$ with $\left\|Y_1\right\|_{\infty} = \pi$ and a unitary eigenvector $\xi_1$ of $Y_1$ corresponding to an eigenvalue $\la_1$ such that $|\la_1|=\pi$, such that,
$$
\alpha(t)=\begin{pmatrix}
e^{it\la_1}&0\\
0&\alpha_1(t)
\end{pmatrix}\begin{array}{l}
\ese_1\\
\ese_1^\bot
\end{array},
$$
where $\ese_1=\mbox{span}\{\xi_1\}$, and $\alpha_1:[0,1]\to \cU(\ese_1^\bot)$ is any curve joining $\sub{I}{\ese_1^\bot}$ and $\sub{U}{\ese_1^\bot}=\left.U\right|_{\ese_1^\bot}$ such that $\|\dot\alpha_1\|_{\infty}\leq \|X\|_{\infty}$. 
 
\medskip

If $\alpha_1(1) = e^{iX_1}$ is such that $\left\|X_1\right\|_{\infty} = \pi$, we apply Lemma \ref{CPR2} to the curve $\alpha_1(t)$: there exists $Y_2$ with $\left\|Y_2\right\|_{\infty} = \pi$ and a unitary eigenvector $\xi_2$ of $Y_2$ corresponding to an eigenvalue $\la_2$ with $|\la_2|=\pi$, such that,
$$
\alpha_1(t)=\begin{pmatrix}
e^{it\la_2}&0\\
0&\alpha_2(t)
\end{pmatrix}\begin{array}{l}
\ese_2\\
\ese_2^\bot
\end{array},
$$
where $\ese_2=\mbox{span}\{\xi_2\}$, and $\alpha_2:[0,1]\to \cU(\ese_2^\bot)$ is any curve joining $\sub{I}{\ese_2^\bot}$ and $\sub{U}{\ese_2^\bot}=\left.U\right|_{\ese_2^\bot}$ such that $\|\dot\alpha_2\|_{\infty}\leq \|X\|_{\infty}$. 

\medskip 

If $\alpha_2(1) = e^{iX_2}$ is such that $\left\|X_2\right\|_{\infty} = \pi$, we continue iterating this Lemma. At some point it will finish because we are in finite dimension. In this way we construct $Y_{\ese}\,$.
\edem

\section{The Grassmannian}

The Grassmannian $\mathcal{G}_n$ is the set of subspaces of $\C^{n}$, which can be identified with the set of orthogonal projections in $\mat $. If we consider in $\mat$ the topology defined by any of all the equivalent norms, the Grassmann space endowed with the inherited topology becomes a compact set. However, it is not connected. Indeed, it is enough to consider the trace 
$\tr$, which is a continuous map defined on the whole space $\mat$, and restricted to $\mathcal{G}_n$ 
takes only positive integer values. 
In particular, this shows that the connected components of $\mathcal{G}_n$ are 
the subsets $\mathcal{G}_{m, n}$ defined as:
$$
\mathcal{G}_{m, n} := \left\{P \in \mathcal{G}_n : \tr(P) = m\right\}.
$$
Each of these components is a submanifold of $M_n(\C)$ [18, p. 129], and connected components are given by the unitary orbit of a given projection $P$ such that $\tr(P) = \rk \, P = m$:
$$
\mathcal{G}_{m, n} = \left\{UPU^{*} : U \in \mathcal{U}(n)\right\}.
$$
The tangent space at a point $P \in \mathcal{G}_{m, n}$ can be identified with the subspace of $P$-codiagonal Hermitian matrices, i.e.
$$
T_P\mathcal{G}_n = \left\{X \in \mathcal{H}(n) : X = PX + XP\right\}.
$$
Denote by $\ese= R(P)$. So each $X \in T_P\mathcal{G}_n$ has a block decomposition
$$
X=\begin{pmatrix}
0&A\\
A^*&0
\end{pmatrix}\begin{array}{l}
\ese\\
\ese^\perp
\end{array} \peso{for some} A \in L(\ese\orto \coma \ese) \ .
$$

In particular note that $T_P\mathcal{G}_n$ has a natural complement $N_P$, which is the space of Hermitian matrices that commute with $P$, that is, the $P$-diagonal Hermitian matrices. The decomposition in diagonal and codiagonal matrices defines a normal bundle, and leads to a covariant derivative
$$
\nabla_V \Gamma(P) = \Pi_{T_P || N_P} \frac{d}{dt} \Gamma(\alpha(t))|_{t=0}, 
$$
where $\Gamma$ is a vector field along the curve $\alpha : (-\varepsilon, \varepsilon) \rightarrow \mathcal{G}_{m, n}$ that satisfies $\alpha(0) = P$ and $\dot{\alpha}(0) = V$. So, we have a notion of parallelism, and the geodesics in this sense are described by the following theorem:

\begin{teo}\rm (See Porta-Recht \cite[(2) and (4)]{PR}, also 
Davis-Kahan \cite{DK} or Halmos \cite{H})
\ 
\begin{itemize}
\item  The unique geodesic at $P$ with direction $X \in T_P\mathcal{G}_n$ is:
$$
\gamma(t) = e^{itX}Pe^{-itX}.
$$
\item  If $P, Q \in \mathcal{G}_n$ are such that $\left\|P - Q\right\|_{\infty} < 1$, there exists a unique $X\in T_P\mathcal{G}_n$ with $\left\|X\right\|_{\infty} < ~\pi/2$ such that 
\beq
\label{el X}
X = PX + XP \in T_P\mathcal{G}_n \peso{and} Q = e^{iX}Pe^{-iX} \ .
\eeq
Therefore they can be join by the geodesic $\gamma(t) = e^{itX}Pe^{-itX}$.

\end{itemize}
\end{teo}

\pausa
{\bf Finsler metrics on the Grassmannian:} For a given symmetric norm
$\|\cdot\|$ in $\mat$, the Grassmann space carries the Finsler structure given by
$$
\left\|X\right\|_P = \left\|X\right\|_{\infty}
\peso{for}X \in T_P\mathcal{G}_n \ .
$$
With this structure, the Grassmann component $\left\{UPU^{*} : U \in \mathcal{U}(n)\right\}$ is isometric (modulo a factor 2) to the orbit  $\left\{US_PU^{*} : U \in \mathcal{U}(n)\right\}$
of the symmetry $S_P = 2P-I$. A straightforward computation shows that, if $X = XP+PX$, then 
$$
e^{iX}S_P = S_Pe^{-iX} \ .$$ 
Therefore,  our previous results about the unitary group can by applied to the Grassmann manifolds. 
As the unitary group acts transitively in these components via $U \cdot P = UPU^{*}$, they are also homogeneous spaces of $\mathcal{U}(n)$. They can be distinguished from other homogeneous submanifolds of $\mathcal{U}(n)$, because the map
$$
P \mapsto S_P = 2P - 1
$$
embeds them in $\mathcal{U}(n)$, and the map $S$ is two times an isometry. The images $S_P$ are symmetries, i.e. matrices that satisfy $S_P^{*} = S_P = S_P^{-1}.$

\begin{nota}\label{el SX} \rm 
Let $P\in \mathcal{G}_{m, n}$ and $X\in T_P\mathcal{G}_n \,$ and 
$A \in L(\ese\orto, \ese)$ such that
\beq\label{modulo}
X=\begin{pmatrix}
0&A\\
A^*&0
\end{pmatrix}\begin{array}{l}
\ese\\
\ese^\perp
\end{array} \ .
\peso{Then}
\left|X\right|=\begin{pmatrix}
\left|A^*\right|&0\\
0&\left|A\right|
\end{pmatrix}\begin{array}{l}
\ese\\
\ese^\perp
\end{array}.
\eeq
\ben
\item
If $\sigma(A) = (\sigma_1 \coma \dots \coma \sigma_k)$ are the (non zero) singular values of $A$, then 
$$
\la(X) = (\sigma_1 \coma \dots \coma \sigma_k \coma 0 \, \uno _{n-2k}\coma -\sigma_k \coma \dots \coma -\sigma_1 ) \in (\R^{n})\da \ ,
$$
and the singular values of $X$ come in pares. 
\item 
We will call $P_X \in \mat $ the orthogonal projection onto the subspace 
\beq\label{omegas}
S_X \igdef \ker(\|X\|I-|X|) \stackrel{\eqref{modulo}}{=} \Omega_1 \oplus \Omega_2  \ , 
\eeq
where 
\begin{align*}
\Omega_1 &= \left\{v\in \ese : \left|A^*\right| v = \left\|A\right\|\, v\right\} \peso{and} 
\Omega_2 = \left\{w \in \ese^\perp: \left|A\right| w = \left\|A\right\|\, w\right\}.
\end{align*}

\item 
We will call $\gamma_X : \R \rightarrow  L(S_X)$ the curve given by 
$$
\gamma_X(t) = P_{X}\,\big( e^{itX}Pe^{-itX} \big)\, P_{X} \Big| _{S_X} \peso {for} t \in \R \ .
$$
By Eq. \eqref{omegas},  $P\, P_X =  P_X \,P $ and $X\,P_X= P_X\,X$. Then also  
$$
P_X \, e^{itX}Pe^{-itX}  =  e^{itX}Pe^{-itX} \,P_X \ ,
$$ 
and 
$\gamma_X(t)\in\cG(S_X)$ for every $t\in \R$ (as well as the compression to $S_X\orto$). 
\EOE
\een
\end{nota}

\pausa
In the following statement we shall use the Notation 
\ref{el SX}.

\begin{teo}
Let $P, Q \in \mathcal{G}_n$ with $\left\|P - Q\right\|_{\infty} < 1$. 
Let $X\in T_P\mathcal{G}_n $ as in Eq. \eqref{el X}. Then $\|X\|_{\infty} < 
\pi/2$ so that $\spec{|X|}\subseteq [0,\pi/2)$. 
%
%
If $\spec{|X|}$ 
has more than one element, then the minimal curves joining $P$ with $Q$ have the following structure:
$$
\delta(t)=\begin{pmatrix}
\gamma_X(t)&0\\
0&\alpha_1(t)
\end{pmatrix}\begin{array}{l}
\ese_X\\
\ese_X^\bot
\end{array},
$$
where $\alpha_1:[0,1]\to \mathcal{G}(\ese_X^\bot)$ is any curve joining 
$P\big |_{\ese_X^\bot}$ and $Q\big | _ {\ese_X^\bot}$ such that $\|\dot\alpha_1\|_{\infty} \leq \|X\|_{\infty}$.
\end{teo}

\bdem
Let's take $\delta(t)$ a minimal curve joining $P$ and $Q$ in $\mathcal{G}_n$ (in this case $L(\delta) = \left\|X\right\|_{\infty}$). Note that if $\gamma(t) = e^{itX}Pe^{-itX}$ then
$$
S_{\gamma(t)} = 2\gamma(t) - 1 = e^{itX}S_Pe^{-itX} = e^{2itX}S_P = S_Pe^{-2itX},
$$
and this curve is the geodesic joining $S_P$ with $S_Q$ and has length $2\left\|X\right\|_{\infty}$. Now consider the curve $\alpha(t) = S_{\delta(t)}$ joining $S_p$ and $S_Q$ in $\mathcal{U}(n)$ and note that it is a minimal curve as
$$
\mbox{L}(\alpha(t)) = \int_{0}^{1} \left\|\dot{\alpha}(t)\right\|_{\infty} = \int_{0}^{1} \left\|2\dot{\delta}(t)\right\|_{\infty} = 2\left\|X\right\|_{\infty}.
$$ 
So from Theorem \ref{structure}, 
%
since $\spec{|X|}\subseteq [0,\pi/2)$ and it has more than one element, then the minimal curves joining $S_P$ with $S_Q$ have the following structure:
$$
\alpha(t)=\begin{pmatrix}
S_Pe^{-2it\sub{X}{\ese}} &0\\
0&\alpha_1(t)
\end{pmatrix}\begin{array}{l}
\ese\\
\ese^\bot
\end{array},
$$
where $\ese=\ker(\|X\|_{\infty}I-|X|)$, $\sub{X}{\ese}=\left.X\right|_{\ese}\,$, and $\alpha_1:[0,1]\to \cU(\ese^\bot)$ is any curve joining 
$S_P\big|_{\ese^\bot}$ and  $S_Q\big|_{\ese^\bot}$  such that $\|\dot\alpha_1\|_{\infty}\leq 2\|X\|_{\infty}\,$.
%
It is easy to prove that actually $\mathcal{S} = \mathcal{S}_X\,$. 

\pausa
As 
the map $\Phi : \cG_n \to \matu$ given by $\Phi(P) = S_P = 2P-I$  is two times an isometry and it  satisfies that $\Phi(P)^{2} = I$ for every $P \in \cG_n\,$, 
 the given structure follows.
\edem

\pausa
As before, a direct consequence of this result is the following corollary about the uniqueness of minimal curves:

\begin{cor}\label{unica2}
Given $P, Q \in \mathcal{G}_n$ with $\left\|P - Q\right\|_{\infty} < 1$, there exists a unique minimal curve between them if and only if $\spec{S_QS_P}\subseteq\{e^{i\theta}, e^{-i\theta}\}$ for some $|\theta|<\pi/2$.
\end{cor}

\begin{rem}
The previous corollary tells us some cases where uniqueness is impossible. As we mention before the singular values of $X$ comes in pares, for every real positive eigenvalues there is a real negative eigenvalue, besides of the possibly zero eigenvalues. If the dimension $n$ is odd then zero is also an eigenvalue and the condition $\spec{S_QS_P}\subseteq\{e^{i\theta}, e^{-i\theta}\}$ couldn't be satisfied. If the dimension $n$ is even a necessary condition for uniqueness is that $\dim(R(P)) = \dim(N(P))$. 
\end{rem}

\section{Minimal curves for the trace norm}

The aim of this section is to study the structure of minimal curves 
joining two matrices, but now where the length is measured 
with respect to the trace norm. 
The spaces that we are interest are the unitary group $\matu$, to complete the study of the previous sections, and the space of  all positive invertible matrices $\matppos$, 
studied by Lim \cite{L1} for the spectral norm. 

\pausa
In order to study properties of minimal curves in this spaces we will focus on the study of the (real vector-)space of Hermitian  matrices $\matsa$, because every positive matrix 
$P \in \matppos$ can be written on the form $P = e^{A}$, where $A \in \mathcal{H}(n)$ and similarly, because every unitary matrix $U \in \mathcal{U}(n)$ can be written on the form $U = e^{iB}$  where $B \in \mathcal{H}(n)$.

\subsection{Study of minimal curves in $\matsa$ for the trace norm}

In this section we shall give a characterization of the minimal curves measured with the trace norm in the space $\matsa$. To achieve this it will be useful to recall a fact that is surely known: which are the minimal curves in $\R^{n}$? 

\begin{pro}

Let $\alpha : [0, 1] \rightarrow \R^{n}$ be a smooth curve joining the vector $O = (0, \ldots, 0)$ with $V = (v_1, \ldots, v_n)$ in $\R^{n}$. The following conditions are equivalent

\begin{enumerate}[(a)]
	\item The curve $\alpha : [0, 1] \rightarrow \R^{n}$ is minimal.
	\item The curves $\alpha_j : [0, 1] \rightarrow \R$ joining $0$ with $v_j$ are minimal.
	\item The curves $\alpha_j : [0, 1] \rightarrow \R$ joining $0$ with $v_j$ are such that $\dot{\alpha}_j(t) > 0$, if $v_j$ is strictly positive, or $\dot{\alpha}_j(t) < 0$, if $v_j$ is strictly negative, for all $t \in [0, 1]$. If $v_j = 0$ then $\alpha_j(t) = 0$ for all $t \in [0, 1]$.
\end{enumerate}

\end{pro}

The next result also should be known and says that the segment joining two Hermitian matrix is shorter than any smooth curve joining them, when  measured with the trace norm. In particular, $$d_1(A, B) = \left\|A - B\right\|_1.$$ We give a simple proof of this fact. 

\begin{lem}
The segments are minimal curves in $\matsa$ with respect to the trace norm.
%
%
%
\end{lem}

\bdem
Since every curve can be approximated by polygonals, it suffices to prove that the segment joining the matrices  $0$ and $A \in \mathcal{H}(n)$ is shorter than any polygonal path joining them. And this fact follows inductively from the following statement. The argument can be seen in this picture: 
\begin{figure}[h!]
\centering
\includegraphics[width=13cm]{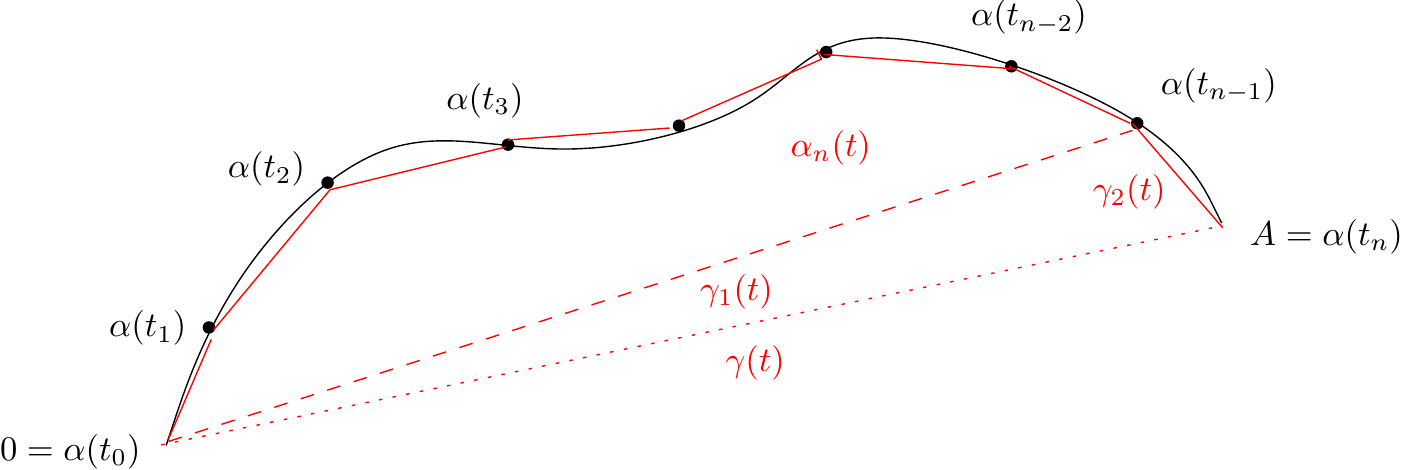} 
\end{figure} 
 \newpage

Let $\gamma(t) = tA$ the segment joining $0$ with $A$ and take any  $B\in \matsa$. 
Then the polygonal  
$$
\gamma_1(t) = \left\{\begin{array}{ccr} 
2tB & \quad & t \in [\,0 \coma \frac12\,] \\  
(2 - 2t)B + (2t - 1)A & \quad &  t\in [\,\frac12 \coma 1\,]
\end{array} \right.
$$
satisfies that 
$$
L(\gamma_1)  = \int_{0}^{1} \left\|\dot{\gamma}_1(t)\right\|_1 dt = \left\|B\right\|_1 + \left\|A - B\right\|_1 \geq \left\|A\right\|_1 = L(\gamma).
$$

\edem

We shall see that a necessary condition for a curve in $\matsa$ 
to be minimal 
is that the entries of the diagonal of the curve should be minimal too. 
In order to state properly this condition 
we recall the definition of the pinching operator. 

\begin{defi}\rm 
Let  $\cP = (P_1\coma \ldots \coma P_k) \in \mat^k$ be a system of projectors in $\mat$. This means that the entries $P_i$ are mutually orthogonal projections 
such that $\sum\limits_{i=1}^k P_i = I.$
The associated pinching operator $\cC_\cP : \mat \to \mat$ is given by 
$$
\mathcal{C}_\cP(A) = \sum_{j = 1}^{k} P_j \, A\,  P_j 
\ , \peso{for every} A\in \mat \ .
$$
\end{defi}

\pausa
One known property of this operator is that, if $\mnorm{\cdot}$ is a 
unitarily invariant norm, then 
\begin{equation}\label{desigualdadpinching}
\mnorm{\mathcal{C}(A)} \leq \mnorm{A} 
\ , \peso{for every} A\in \mat \ .
\end{equation}
Actually, with respect to the trace norm,  we will prove a stronger result which in some sense asserts that the pinching operator preserves minimality of curves:

\begin{lem}\label{pinchingcorta}

Let $X : [0, 1] \rightarrow \mathcal{H}(n)$ be a smooth  curve which is minimal for the trace norm $\left\|\cdot\right\|_1\,$. Assume further that 
$X(0) = 0$. Let  $D = X(1) \in \mathcal{H}(n)$ 
and $\cP$ a system of projectors, all of them commuting with $D$. 
%
Then the pinched curve $\mathcal{C}_\cP (X(t))$ is also a minimal
 curve joining $0$ with $D \in \mathcal{H}(n)$.

\end{lem}

\bdem

It's an immediate consequence of \eqref{desigualdadpinching} and the fact that 
\beq\label{deripinch}
\frac{d}{dt}\mathcal{C}_\cP(X(t)) = \mathcal{C}_\cP(\dot{X}(t)) 
\stackrel{\eqref{desigualdadpinching}}{\implies }
\left\|\mathcal{C}_\cP(\dot{X}(t))\right\|_1 \leq \left\|\dot{X}(t)\right\|_1 \ 
\peso{for} t\in [0\coma 1] \ .
\eeq
Therefore $L(\mathcal{C}_\cP(X(t))) \leq L((X(t)) = \left\|D\right\|_1\,$. On the other hand, $\mathcal{C}_\cP(X(1)) = \mathcal{C}_\cP(D) = D$ by the commutativity hypothesis.
\edem

A particular case is the following: Let $\cP_0$ be the system of projectors
associated to the canonical basis of $\C^n$, so that the image of $\cC_{\cP_0}$ 
is the set of usual diagonal matrices.   

\begin{cor}\label{diagonal} \rm

Let $X : [0, 1] \rightarrow \mathcal{H}(n)$ be curve as in Lemma \ref{pinchingcorta}.
Assume further that $D = X(1) = \diag {\, d_1 \coma \dots \coma d_n \, } $ is a diagonal matrix.  Then, for every $t\in  [0\coma 1] $, 
$$ 
\cC_{\cP_0}(X(t)\,) = \diag {\, x_{1}(t) d_1 \coma \dots \coma x_{n}(t) d_n \, }
,
$$
where $x_i : [0, 1] \rightarrow \R$ is such that $x_i(0) = 0, x_i(1) = 1$ and $\dot{x}_i(t) > 0$ if $d_i$ is positive, or  $\dot{x}_i(t) < 0$ if $d_i$ is negative and $x_i(t) = 0$ for all $t \in [0, 1]$ if $d_i = 0$. In other words, the curves in the diagonal of the matrix are also minimal.
\end{cor}

We are almost ready for the characterization of the minimal curves in $\matsa$ but, before that, we prove a lemma on which is based mainly the result. 

\begin{lem}\label{dospordos}

Let $X : [0, 1] \rightarrow \mathcal{H}(2)$ be a (smooth) minimal curve, measure with the trace norm $\left\|\cdot\right\|_1$, joining $0$ with $D \in \mathcal{H}(2)$,
$$
D = \begin{pmatrix}
\alpha & 0 \\ 0 & \beta
\end{pmatrix}.
$$
\begin{enumerate}[(a)]
	\item If $\alpha, \beta \geq 0$, then $\dot{X}(t) \geq 0$.
	\item If $\alpha, \beta < 0$, then $\dot{X}(t) \leq 0$.
	\item If $\alpha > 0, \beta < 0$, then $X(t)$ is a diagonal matrix.
\end{enumerate}

\end{lem}

\proof
\begin{enumerate}[(a)]
\item Note that, by Lemma \ref{pinchingcorta}, the diagonal curve 
$$
\mathcal{C}(X(t)) = \begin{pmatrix} x_{11}(t) & 0 \\ 0 & x_{22}(t) \end{pmatrix}
$$
is  minimal. Using Eq. \eqref{deripinch}, we get that 
\begin{equation}\label{eq01}
\tr  \ |\mathcal{C}(\dot{X}(t))| = \tr  \ |\dot{X}(t)|,
\end{equation}
as the curve is smooth and $L(\mathcal{C}(X(t))) = L(X(t))$. On the other hand, 
by Corollary \ref{diagonal} we know 
that  $\dot{x}_{11}(t)\ge 0 $ and $ \dot{x}_{22}(t) \geq 0 \implies\mathcal{C}(\dot{X}(t)) \geq 0$. If we combine
$$
\tr \ |\mathcal{C}(\dot{X}(t))| = 
\tr \ \mathcal{C}(\dot{X}(t)) = 
\tr \ \dot{X}(t) \leq 
\tr \ |\dot{X}(t)|,
$$
with \eqref{eq01} we have $\tr \ \dot{X}(t) = \tr \ |\dot{X}(t)|$, which implies that $\dot{X}(t) = \left|\dot{X}(t)\right| \geq 0$.

\item It follows from the previous case taking $Y(t) = -X(t)$.

\item Let
$$
X(t) = \begin{pmatrix} a(t) & b(t) \\ \overline{b(t)} & c(t) \end{pmatrix},
$$
and let $\lambda_1(t), \lambda_2(t)$ be the eigenvalues of $\dot{X}(t)$ arranged in decreasing order. As
$$
(\dot{a}(t), \dot{c}(t)) \prec (\lambda_1(t), \lambda_2(t)),
$$
and $\dot{a}(t) \geq 0, \dot{c}(t) \leq 0 $ (by Corollary \ref{diagonal}) then
\begin{equation}\label{mayo}
\lambda_2(t) \leq \dot{c}(t) \leq 0 \leq \dot{a}(t) \leq \lambda_1(t)
\peso{and} \left|\dot{c}(t)\right| + \left|\dot{a}(t)\right| \le \left|\lambda_1(t)\right| + \left|\lambda_2(t)\right|  
\end{equation}
for every $t\in [0\coma 1]$. 
On the other hand, using Lemma \ref{pinchingcorta} we get that 
\begin{align*}
\int_{0}^{1} \left|\dot{c}(t)\right| + \left|\dot{a}(t)\right|dt & =	\int_{0}^{1} tr(|\mathcal{C}(\dot{X}(t))|)dt = \int_{0}^{1}  tr(|\dot{X}(t)|)dt = \\ & = \int_{0}^{1} \left\|\dot{X}(t)\right\|_1 dt = \int_{0}^{1} \left|\lambda_1(t)\right| + \left|\lambda_2(t)\right|dt,
\end{align*}
As the curve is smooth, by Eq. \eqref{mayo}  
$$
\left|\dot{c}(t)\right| + \left|\dot{a}(t)\right| = \left|\lambda_1(t)\right| + \left|\lambda_2(t)\right|  \implies 
\dot{a}(t) = \lambda_1(t), \ \ \ \ \dot{c}(t) = \lambda_2(t).
$$
Applying the Frobenius norm to the matrix $\dot{X}(t)$ 
we conclude that $\dot{b}(t) = 0$ for every $t\in [0\coma 1]$. 
Since $b(1) = 0$ we conclude that  $b(t) = 0$ for every $t\in [0\coma 1]$.
\QED
\end{enumerate}

\pausa
All minimal curves should have the following structure: 

\begin{teo}\label{caracterizacion}
Let $X : [0, 1] \rightarrow \mathcal{H}(n)$ be a smooth  curve such that 
$X(0) = 0$ and $D = X(1) = \diag {\, d_1 \coma \dots \coma d_n \, } $ is a diagonal matrix such that 
%
%
its diagonal entries $d_j$ are arranged in decreasing order. 
Then, $X$ is a minimal curve joining $0$ with $D$ with respect to the trace norm  if and only if it has the form:
\ben
\item If 
$\ese_1 = span \left\{e_j :  d_j > 0\right\}  \coma 
\ese_2  = span \left\{e_j : d_j = 0\right\}$ and
$\ese_3  = span \left\{e_j : d_j < 0\right\}$, then 
$$
X(t)=\begin{pmatrix}
P(t)&0&0\\
0&0&0\\
0&0&N(t)
\end{pmatrix}\begin{array}{l}
\ese_1\\
\ese_2\\
\ese_3
\end{array},
$$
\item  $P : [0, 1] \rightarrow \mathcal{H}(\ese_1)$ is a minimal curve measure with the trace norm joining $\sub{0}{\ese_1}$ with $\sub{D}{\ese_1}$ such that $\dot{P}(t) \geq 0$ for all $t \in [0, 1]$.  
\item Also $N : [0, 1] \rightarrow \mathcal{H}(\ese_3)$ is a minimal curve measure with the trace norm joining $\sub{0}{\ese_3}$ and  $\sub{D}{\ese_3}$ such that $\dot{N}(t) \leq 0$ for all $t \in [0, 1]$. 
\een
%
\end{teo}
\proof
 It's easy to check that if $X(t)$ has the given form then 
\begin{align*}
L(X) &= \int_{0}^{1} \left\|\dot{X}(t)\right\|_{1}dt = \int_{0}^{1} \left\|\dot{N}(t)\right\|_{1}dt + \int_{0}^{1} \left\|\dot{P}(t)\right\|_{1}dt \\ 
& = \tr \int_{0}^{1} \dot{P}(t) -\dot{N}(t) \, dt = \tr \big(P(1)-N(1)\,\big)
= \|D\|_1 
= L(t\, D),
\end{align*}
so it's a minimal curve. 
To see the converse the key is the previous Lemma \ref{dospordos}. Let 
$$
X(t)=\begin{pmatrix}
X_{11}(t)&X_{12}(t)&X_{13}(t)\\
X_{21}(t)&X_{22}(t)&X_{23}(t)\\
X_{31}(t)&X_{32}(t)&X_{33}(t)
\end{pmatrix}\begin{array}{l}
\ese_1\\
\ese_2\\
\ese_3
\end{array},
$$
then
\begin{enumerate}[(a)]

	\item \textit{Claim:} $\dot{X}_{11}(t) \geq 0, X_{22}(t) = 0$ and $\dot{X}_{33}(t) \leq 0$. Indeed if $$\mathcal{C}(X(t)) = \sum_{j = 1}^{3} P_j X(t) P_j,$$ where $P_j$ is the projection over $\mathcal{S}_j$, by Lemma \ref{pinchingcorta}, the block $X_{22}(t) = 0$ (because it has to be minimal joining $0$ with $0$). A consequence of Lemma \ref{dospordos}(a) is that $\dot{X}_{11}(t) \geq 0$ because is a minimal curve joining $\sub{0}{\ese_1}$ with $\sub{D}{\ese_1}$ (the positive entries). It's natural to call it $X_{11}(t) = P(t)$. Analogously  $\dot{X}_{33}(t) \leq 0$ so it's natural to call it $X_{33}(t) = N(t)$.
	
	\item \textit{Claim:} $X_{12}(t) = X_{21}(t) = 0$. Let's take an appropriate pinching operator in this way: choose a diagonal entry
	of the block matrix $P(t)$ and call it $p_{11}(t)$. This entry is associate to an eigenvalue $d_i > 0$. Now choose a diagonal entry
of the block matrix $0$. 
This entry is associate to an eigenvalue $d_j = 0$. Let 
$\cM_{i} = span \left\{e_i \right\}$,  
$\cM_{j} = span \left\{e_j \right\}$ and 
$P$ be the projection onto $\cM_{i} \oplus \cM_{j}$. 
This projection together with the orthogonal projection onto 
$(\cM_{i} \oplus \cM_{j})\orto$ 
will produce the desired pinching operator (which commutes with $D$
because it is diagonal). 
Using a permutation matrix we can have in the first block the following matrix:
$$
Y(t) =\begin{pmatrix}
p_{11}(t)&y_{12}(t)\\
y_{21}(t)&0
\end{pmatrix}\begin{array}{l}
\eme_i\\
\eme_j
\end{array}, 
$$
which is minimal joining $0$ and $\diag {d_i\coma 0}$ in $\cH(2)$, by Lemma \ref{pinchingcorta}. But $\dot{Y}(t) \geq 0$ by Lemma \ref{dospordos}(a), so  $\dot{y}_{12}(t) = \dot{y}_{21}(t) = 0$ and therefore $y_{12}(t) = y_{21}(t) = 0$. Using this trick we can put a zero in every entry of $X_{12}(t)$ and every entry of $X_{21}(t)$.

	\item \textit{Claim:} $X_{23}(t) = X_{32}(t) = 0$. It's deduce from the previous case as before, take
$$
Y(t)= -\begin{pmatrix} 0 & y_{23}(t) \\ y_{32}(t) & n_{22}(t) \end{pmatrix}\begin{array}{l}
\eme_i\\
\eme_j
\end{array}.
$$

\item \textit{Claim:} $X_{13}(t) = X_{31}(t) = 0$. 
It follows by a restriction similar as before, using in this case 
Lemma \ref{dospordos}(c). Choosing $e_i \in \ese_1$ and $e_j\in \ese_3$ 
one can get a block of the form 
$$
Y(t) =\begin{pmatrix}
p(t)&y_{12}(t)\\
y_{21}(t)&n(t)
\end{pmatrix}\begin{array}{l}
\eme_i\\
\eme_j
\end{array},  \peso{with} p(t) >0 \peso{and} n(t)<0 \ ,
$$
which is minimal by Lemma \ref{pinchingcorta}. In this case  $y_{12}(t) = y_{21}(t) = 0$ by Lemma \ref{dospordos}(c). In this way we find a zero in every entry of $X_{13}(t)$ and every entry of $X_{31}(t)$.
\QED
\end{enumerate}

\begin{rem}
Using the hypothesis and notations of Theorem \ref{caracterizacion}, 
since $P(0) = 0$ and $\dot{P}(t) \geq 0$ for every $t\in [0\coma 1]$, then 
the curve $P(t)$ is increasing in $\matsa$, and in particular every $P(t) \ge 0$.
Similarly, the curve $N(t) $ is decreasing and $N(t)\le 0$. 
%
%
\EOE
\end{rem}

A consequence of this characterization is the following necessary condition 
about the curves of eigenvalues. 

\begin{cor}\label{eigenvalues hermitian}
Let $X : [0, 1] \rightarrow \mathcal{H}(n)$ be a (smooth) minimal curve, measure with the trace norm, joining $0$ with $A\in \matsa$. 
Assume that the (continuous) curves $\la_j(X(t))$ joining $0$ with $\lambda_j(A)$ are piecewise smooth. Then they are minimal curves (in $\R$)  
with respect to the trace norm. By Theorem 
\ref{caracterizacion} this means that they are monotone maps.  

\end{cor}

\bdem
By a unitary conjugation (which preserve all eigenvalues), we can assume that 
$A =D= \diag {\, d_1 \coma \dots \coma d_n \, }$, a diagonal matrix as 
in Theorem \ref{caracterizacion}. 
As $X(t)$ is a minimal curve it must have the form given in the previous Theorem \ref{caracterizacion}, so the positive eigenvalues correspond to $P(t)$ and the negatives to $N(t)$. Suppose that $\dim(\ese_1) = k$ and fix $j\le k$.

\pausa
Note that $L(\lambda_j(P(t))) \geq d_j$ because the curve joins $0$ with $d_j$ 
and $L(t\,d_j) = 
L(p_{jj}(t)) = d_j\,$, which are minimal by Corollary \ref{diagonal}. 
On the other hand, as $d(P(t)) \prec \lambda(P(t))$, 
$$
\sum_{j = 1}^{k} p_{jj}(t) = \sum_{j = 1}^{k} \lambda_j(P(t)).
$$
If we take derivative and then integrate, we have
$$
\sum_{j = 1}^{k} d_j  = \sum_{j = 1}^{k} L(\lambda_j(P(t))),
$$
i.e. $L(\lambda_j(P(t))) = d_j$, so they are minimal curves. With minor changes it can be prove that $L(\lambda_j(N(t))) = -d_j$ if $e_j \in \ese_3\,$.
\edem

\subsection{The space $\matinv^+ $ of positive invertible matrices}

\medskip
Recall that  $\matppos$ is an open subset of the space of hermitian matrices $\matsa$. Therefore, it inheres a geometric structure where the tangent spaces can be identified with $\matsa$. Also, there exists a natural transitive action of $\matinv$ by conjugation. The properties of this action make $\matppos$ become an homogeneous space. Moreover, using this action it is possible to define a covariant derivative that leads to the following differential equation for the geodesics
$$
\gamma''=\gamma'\gamma^{-1}\gamma',
$$
 (see \cite{CoMa}, \cite{CPR3} and \cite{KN}). Given $A,B\in\matppos$, the solution of the corresponding Dirichlet problem gives the following expression for the geodesic joining $A$ with $B$
\beq\label{geos}
\sub{\gamma}{\matppos}(t)=A^{1/2}(A^{-1/2}BA^{-1/2})^tA^{1/2}.
 \eeq
Note that if one of the end points is the identity, also in this case the geodesic is a one-parameter subgroup of $\matinv$. 

\medskip
It is also possible to define  a (canonical) Riemannian structure on $\matppos$. Indeed, consider the inner product associated to the trace in the tangent space of $\matppos$ at the identity. Then, using the homogeneous structure, we can define the following inner product
$$
\pint{X,Y}_A=\tr\big((A^{-1/2}XA^{-1/2})(A^{-1/2}YA^{-1/2})^*\big)=\tr(A^{-1/2}XA^{-1}Y^*A^{-1/2}),
$$
in the tangent space corresponding to another  point $A\in\matppos$. Endowed with this structure, the action by conjugations becomes isometric, and $\matppos$ becomes a complete Riemannian manifold with non-positive sectional curvature. As before, these geodesics are minimal curves. However, in this case, they are short not only for $t\in [0,1]$, but also for every $t\in\R$  (see \cite{Bh PM}, \cite{PR}).

\pausa
Simlarly, we have defined the invariant Finsler metric associated to the trace norm, by 
$$
\sub{\|X\|}{1}=\|A^{-1/2}XA^{-1/2}\|_1 = \tr\, |A^{-1/2}XA^{-1/2}| \ ,
$$
for any $X\in \matsa$ thought as a tangent vector at $A\in \matppos$. 
So, given an interval $[a, b] \subset \R$ and a smooth curve $\alpha: [a, b] \rightarrow \matppos$, the length of this curve is defined as
$$
L(\alpha) = \int_{a}^{b} 
\left\|\alpha^{-1/2}(t)\dot{\alpha}(t)\alpha^{-1/2}(t)\right\|_{1}dt \ .
$$ 
Here the geodesics given in Eq. \eqref{geos} are still minimal, but they are not unique.  
Remember that the space $\matinv^+$ of positive invertible matrices is a nonpositive 
curvature space so the exponential map increases distances in general. 
Hence it's natural to think that the logarithm of a minimal curve 
in $\matinv^+ $ could be  a minimal curve in $\matsa$. What is more surprising is that, for the metric induced by the trace norm,  
if we take a minimal curve in $\matsa$, then its exponential gives a minimal curve in $\matinv^+ $.

\medspace

Let's recall the following formula for the derivative (attributed  to Duhamel, Dyson, Feynman, and
Schwingerof) of the exponential map:
\beq\label{der exp}
De^{X}(Y)= \int_{0}^{1} e^{tX} Y e^{(1-t)X}dt \ , \peso{for} X\coma Y \in \mat \ .
\eeq
We also recall a simple proof of the generalized \textit{infinitesimal exponential metric increasing property} (IEMI). In order to prove this result we will use the following theorem, which is stated, for example, in Bhatia's book \cite[Thm. 5.4.7]{Bh PM}.  
\begin{teo} \rm
Let $X\in \mat$ and $A\coma B \in \matpos$. Then, for every unitarily invariant norm, 
$$
\mnorm{A^{1/2}XB^{1/2}}  \leq \mnorm{\int_{0}^{1} A^{t}XB^{1-t}dt}.
$$
\end{teo}

\begin{pro}
(Generalized IEMI) Given $H \coma K$ in $\mathcal{H}(n)$, they satisfies that
$$
\mnorm{e^{-H/2}De^{H}(K)e^{-H/2}} \geq \mnorm{K} 
$$
for every unitarily invariant norm.
\end{pro}

\bdem

As $K = e^{H/2}(e^{-H/2}Ke^{-H/2})e^{H/2}$, then
\begin{align*}
\mnorm{K} & = \mnorm{e^{H/2}(e^{-H/2}Ke^{-H/2})e^{H/2}} \leq \mnorm{\int_{0}^{1} e^{tH}(e^{-H/2}Ke^{-H/2})e^{(1-t)H}dt} = \\
& = \mnorm{e^{-H/2} \left(\int_{0}^{1} e^{tH}Ke^{(1-t)H} dt\right)e^{-H/2}} = \mnorm{e^{-H/2} De^{H}(K) e^{-H/2}}.
\end{align*}

\edem

An immediate consequence of the generalized EMI for the particular case of the trace norm implies that:

\begin{cor}
\label{cor pabajo}
Let $\alpha : [0, 1] \rightarrow \matppos $ be a smooth curve 
such that $\alpha (0) =I$. 
Let $X: [0, 1] \rightarrow \mathcal{H}(n)$ be smooth and such that  
$\alpha(t) = e^{X(t)}$ for $t\in [0\coma 1]$ and $X(0)=0$. 
Then 
\ben
\item We have that 
$L(\alpha)$ (in $\matppos $)  $\ge L(X)$ (in $\matsa$).
\item If $\alpha(t)$ were minimal, then $X(t)$ is also minimal, 
with respect to any fixed UIN.  
\een
\end{cor}

\bdem
By the chain rule $\dot{\alpha}(t) = De^{X(t)}(\dot{X}(t))$. So that
$$
L(\alpha) = \int_{0}^{1} \mnorm {\, e^{-X(t)/2}(t)\, 
De^{X(t)}(\dot{X}(t))\, e^{-X(t)/2}(t)\,} \ dt
\geq \int_{0}^{1} \mnorm{\, \dot{X}(t)\,} dt = L(X) 
$$
by the Generalized IEMI. If $\alpha$ were minimal, let $A = \alpha(1)$ and 
$Y= \log \, A\in \matsa$. Hence $A =e^Y$ and $\mnorm{Y} = L(e^{t\, Y})=L(\alpha)$, since both curves are minimal joining the same points.  
Note that $X(1)  = Y$ (in $\matppos$ the $\log$ is unique), so that 
$$
\mnorm{Y} = L(e^{t\, Y})=L(\alpha) \ge L(X)\ge  L(tY) = \mnorm{Y} \ , $$ 
since $t\, Y$ is minimal joining $0$ and $Y$. This fact also show that $X$ must be minimal. 
\edem

Surprisingly, when restricted to the trace norm, a converse result also holds:
if $X(t)$ is a minimal curve then $\alpha(t)= e^{X(t)}$ is also minimal. 
The proof is based on the characterization given in the previous section. 
We recall the notations $P(t)$, $N(t)$, and $\ese_i$ from Theorem \ref{caracterizacion}, 
which we shall use in the following proofs. We need also to recall some technical results. 

\begin{rem}\label{[1]}
Let $I \subset \R$ be an open interval and let $C^{1}(I)$ be the space of continuously differentiable real functions on $I$. 
Denote by $\matsai = \{A\in \matsa: \sigma (A) \subseteq I\}$. 
Fix $f \in C^{1}(I) $  and consider the smooth map $f : \matsai \to \matsa$, acting  
by functional calculus. 

\pausa
Given $A \in \matsai$,  the derivative 
$Df(A)$ of $f$ at $A$ is a linear map from $\matsa$ into itself. 
If $A= \diag {\, \la_1 \coma \dots \coma \la_n \, }$ a diagonal matrix, there is a formula for $Df(A)$:
\begin{equation}\label{formulita}
Df(A)(B) = f^{[1]}(A) \circ B \ , \peso{for every} B\in \matsa \ ,
\end{equation}
where $\circ$ denotes the Hadamard product and 
$f^{[1]}(A)\in \matsa $ is 
defined as
\begin{align*}
f^{[1]}(A)_{ij} & = \frac{f(\la_i) - f(\la_j)}{\la_i - \la_j} 
& \peso{if} \la_i  \neq \la_j \ \ \\
f^{[1]}(A)_{ii}& = \quad f^{'}(\la_i) & \ \ \ \peso{if} \la_i  = \la_j \ .
\end{align*}
This is called the first divided difference of $f$ at $A$. It is well known 
(see \cite{bhatiapositivas}) that if $f = \exp$, 
then $ f^{[1]}(A) \in \matpos$ for every (diagonal) $ A\in \matsa$. \EOE
\end{rem}

\begin{teo}\label{exp corta}
%
Let $\alpha(t) = e^{Y(t)}$ where $Y: [0, 1] \rightarrow \mathcal{H}(n)$. Then, 
if $Y(t)$ is a minimal curve then $\alpha(t)$ is minimal in $\matppos$.

\end{teo}

\bdem
Let $X: [0, 1] \rightarrow \matsa$ be  any smooth curve. 
By the chain rule $\dot{\alpha}(t) = De^{X(t)}(\dot{X}(t))$.
For any fixed $t\in [0, 1] $,  choose an orthonormal basis such that
$X(t) = \diag {\lambda_1(t) \coma \dots \coma \lambda_n(t)}$. 
By the formula \eqref{formulita}, in this basis, 
$$
De^{X(t)}(\dot{X}(t)) = \exp^{[1]}(X(t)\,) \circ  \dot{X}(t) \ .
$$
In particular, its diagonal entries are $e^{\la_i(t)} \, \dot{X}_{ii}(t)$.  Therefore, the diagonal entries  of 
$$
M \igdef  \alpha^{-1/2}(t)\dot{\alpha}(t)\alpha^{-1/2}(t)
=e^{-X(t)/2} De^{X(t)}(\dot{X}(t))e^{-X(t)/2}  
$$ 
are $M_{ii} = e^{-\la_i(t)/2}( e^{\la_i(t)} \, \dot{X}_{ii}(t)\, ) e^{-\la_i(t)/2} 
= \dot{X}_{ii}(t)$. 
In particular, 
\beq\label{tr L}
\tr \, \big( e^{-X(t)/2} De^{X(t)}(\dot{X}(t))e^{-X(t)/2}  \big)  = \sum_{j = 1}^{n} \dot{X}_{jj}(t) = \tr \, \dot{X}(t) \ ,
\eeq
(and we can forget the change of basis at any $t$). 
Consider now the diagonal decomposition  of $Y(t)$ of Theorem \ref{caracterizacion} 
in terms of $P(t)$ and $N(t)$. Note that, by Remark \ref{[1]},  
$$
\dot{P}(t) \ , \ \ \exp^{[1]}(P(t)) \ , \ \ \exp^{[1]}(N(t)) 
\peso{and} -\dot{N}(t) \in \matpos  
$$ 
(changing basis to get diagonal representations if necessary). Then, 
by Hadamard theorem 
\beq\label{PN}
De^{N(t)}(\dot{P}(t)) \ge 0 \peso{and} De^{N(t)}(\dot{N}(t)) \leq 0  \ .  
\eeq
Finally
\begin{align*}
L(\alpha) & = \int_{0}^{1} \left\|\alpha^{-1/2}(t)\dot{\alpha}(t)\alpha^{-1/2}(t)\right\|_1dt  \\
& = \int_{0}^{1} \left\|e^{-Y(t)/2} De^{Y(t)}(\dot{Y}(t))e^{-Y(t)/2}\right\|_1dt  
=  \int_{0}^{1} \tr \left|e^{-Y(t)/2} De^{Y(t)}(\dot{Y}(t))e^{-Y(t)/2}\right|dt  \\
& =  \int_{0}^{1} \tr \left|e^{-P(t)/2} De^{P(t)}(\dot{P}(t))e^{-P(t)/2}\right|dt 
+ \tr \left|e^{-N(t)/2} De^{N(t)}(\dot{N}(t))e^{-N(t)/2}\right|dt  \\
& \stackrel{\eqref{PN}}{=} 
 \int_{0}^{1} \tr \left(e^{-P(t)/2} De^{P(t)}(\dot{P}(t))e^{-P(t)/2}\right)dt 
- \tr \left(e^{-N(t)/2} De^{N(t)}(\dot{N}(t))e^{-N(t)/2}\right)dt  \\
& \stackrel{\eqref{tr L}}{=} 
 \int_{0}^{1}  \tr \, \dot{P}(t)  - \tr \, \dot{N}(t) \, dt 
= \int_{0}^{1}  \tr \, \left|\dot{Y}(t)\right|dt 
=L(Y) \ . 
\end{align*}
This equality jointly with item 1. of Corollary \ref{cor pabajo} 
easily implies the minimality of $\alpha$. 
\edem
\begin{cor}\label{equiv traza}
Let $\alpha : [0, 1] \rightarrow \matppos$ be a (smooth) minimal curve
joining $I$ with $D\in \matppos$. Denote by $Y(t) = \log \alpha (t) \in \matsa $ 
for $t \in [0, 1]$. Then $\alpha $ is a minimal curve for the trace norm 
$\iff$ $Y$ is minimal in $\matsa$ $\iff$ $Y$ has the block diagonal form given 
in Theorem \ref{caracterizacion} (modulo a change of basis which diagonalizes $D$).
\end{cor}
\proof
It is a consequence ot Corollary \ref{cor pabajo} and Theorem \ref{exp corta} 
\QED

\pausa 
Another necessary condition for a curve to be minimal in this space is that the curves that joins the eigenvalues are also minimal.

\begin{cor}\label{eigenvalues positive}
Let $\alpha : [0, 1] \rightarrow \matppos$ be a (smooth) minimal curve, measure with the trace norm, joining $I$ with $P\in \matppos$. 
Assume that the (continuous) curves $\la_j(\alpha (t))$ joining $0$ with $\lambda_j(P)$ are piecewise smooth. Then they 
are minimal curves (in $\R$)  with respect to the trace norm.
By Theorem 
\ref{caracterizacion} this means that they are monotone increasing.  
%
%
%
%
\end{cor}

\bdem
Let $X: [0, 1] \rightarrow \mathcal{H}(n)$ be the smooth curve given by 
$X(t) = \log \, \alpha(t)$ (so that  $\alpha(t) = e^{X(t)}$) for $t\in [0\coma 1]$. 
Note that $\alpha(t)$ is  minimal  $\iff X(t)$ is minimal. If $X(t)$ is minimal
then the curves  $\lambda_j(X(t)) = \log \, \la_j(\alpha (t))$ are also 
piecewise smooth and therefore minimal, by Corollary \ref{eigenvalues hermitian}. 
Thus the curves 
$e^{\lambda_j(X(t))} = \la_j(\alpha(t))$ are minimal (by the $n=1$ case
ot Corollary \ref{equiv traza}).
\edem

\subsection{The unitary group $\matu$.}
The space $\matu$  has positive 
curvature, so that the exponential map decreases distances in general. 
Hence it's natural to think that the exponential of a minimal curve 
in $i \matsa$ could be  a minimal curve in $\matu$. 
%

\pausa
When considering the trace norm both in $\matsa$ and in $\matu$, we can use the characterization of minimal curves in $\matsa$ to construct 
several examples of minimal curves in the space $\matu$. 
In other words, we have a sufficient condition for a curve 
in $\matu$ to be minimal for the trace norm: Being the exponential 
of a minimal curve in $i\matsa$. 

\begin{teo}\label{corta Un}
Let $X: [0, 1] \rightarrow \mathcal{H}(n)$ be a smooth curve. Let 
$\beta: [0, 1] \rightarrow \matu$ be the smooth curve given by 
$\beta(t) = e^{i X(t)}$  for $t \in [0, 1]$. Then
\ben
\item $L(\beta)  \le L(X)$  (both measured with the trace norm). 
\item If $Y: [0, 1] \rightarrow \mathcal{H}(n)$ is smooth, minimal in $\matsa$, 
$Y(0)=0$ and $\left\|Y(1)\right\|_{\infty} \le \pi$, then 
its exponential $\alpha(t)
 = e^{i Y(t)}$ is a minimal curve in $\matu$.
\een
\end{teo}

\bdem
The first inequality follows because 
\begin{align*}
L(\beta) 
& = \int_{0}^{1} \left\|\dot{\beta}(t)\right\|_1 dt  
\stackrel{\eqref{der exp}}{=}
\int_{0}^{1} \left\|\int_{0}^{1} e^{isX(t)}\dot{X}(t)e^{i(1-s)X(t)}ds \right\|_1 dt \\ 
& \leq  \int_{0}^{1} \int_{0}^{1} \left\|e^{isX(t)}\dot{X}(t)e^{i(1-s)X(t)}\right\|_1 ds dt =  	\int_{0}^{1} \left\|\dot{X}(t)\right\|_1  dt =
	L(X) \ .
\end{align*}
Denote by $Z =Y(1)$ and $U = e^{itZ} = \alpha(1)$. 
Since $\|Z\|_{\infty}\le \pi$,  the curve $\gamma (t) = e^{itZ}$ for $t \in [0, 1]$ 
is  minimal in $\matu$,  joining $I=\alpha(0)$ with $ U= \alpha(1)$, and 
$L(\gamma) = \|Z\|_1\,$. 

\pausa
On the other hand, the curve $t\mapsto tZ$ is minimal in $\matsa$, and it joins 
$0=Y(0)$ with $Z=Y(1)$. Hence $\|Z\|_1 = L(tZ) = L(Y)$, because $Y$ also is minimal.
Then, by item 1,  
$$
L(\alpha) \le L(Y) = \|Z\|_1 = L(\gamma) \ . 
$$
This proves that 
$\alpha(t)$ is another minimal curve joining $I$ with $ U$. 
\edem

On the other hand, a necessary condition for a curve in $\matu$ to be minimal  for the trace norm is that all its curves of eigenvalues are minimal curves
in $\T =\cU(1)$.

\begin{pro}\label{eigenvalues unitary}

Let $\alpha : [0, 1] \rightarrow \mathcal{U}(n)$ be a (smooth) minimal 
curve for the trace norm, 
joining $I$ with $U = e^{iX}$, where $\left\|X\right\|_{\infty} < \pi$. 
Let $Y: [0, 1] \rightarrow \matsa$ be the unique smooth curve such that
$\alpha(t) = e^{iY(t)}$ and $\|Y(t)\|_{\infty} < \pi$ for every $t \in [0, 1]$. 
Assume that the  curves $\la_j(\alpha (t)\,) \igdef 
e^{i\, \la_j(Y(t)\,)}$ for $t \in [0 \coma 1]$ 
are piecewise smooth.   
Then 
\ben
\item The curves $\la_j(Y(t)\,)$ are monotone, so they are minimals in $\R$ 
for the trace norm. 
\item The curves  $\la_j(\alpha (t))$ are minimal in $\C$ for 
trace norm.
\een
\end{pro}
\bdem

%

Fix $t \in [0, 1)$ and denote  
$Y_t = Y(t)$. Since the curve $\alpha$ is minimal and $\|Y_t\|_{\infty}<\pi$, 
\begin{align*}
\left\|X\right\|_1 & = d_1(I, e^{iX}) = d_1(I, e^{iY_{t}}) + d_1(e^{iY_{t}}, e^{iX}) 
\\
& = \left\|Y_t\right\|_1 + d_1(I, e^{-iY_{t}}e^{iX}) 
= \left\|Y_t\right\|_1 + d_1\left(I, e^{i(U_tXU_t^{*} - V_tY_tV_t^{*})}\, \right) \\
& =  \left\|Y_t\right\|_1 + \left\|U_tXU_t^{*} - V_tY_tV_t^{*}\right\|_1 \\
&\geq \left\|Y_t\right\|_1 + \displaystyle\sum_{j = 1}^{n} \left|\lambda_j^{\downarrow}(U_tXU_t^{*}) - \lambda_j^{\downarrow}(V_tY_tV_t^{*})\right|  \\
& =  \left\|Y_t\right\|_1 + \displaystyle\sum_{j = 1}^{n} \left|\lambda_j^{\downarrow}(X) - \lambda_j^{\downarrow}(Y_t)\right| \\
&\geq 
\left\|Y_t\right\|_1 + \displaystyle\sum_{j = 1}^{n} \left|\lambda_j^{\downarrow}(X)\right| - \left|\lambda_j^{\downarrow}(Y_t)\right| = \left\|X\right\|_1.
\end{align*}
For the equality between the distance and the trace norm see for example \cite{ALV}. Besides the first inequality follows from Lidskii's Theorem and the convexity of taking moduli 
(see \cite[III.4.4]{Bh MA}). In conclusion, 
$$
\left|\lambda_j^{\downarrow}(X) - \lambda_j^{\downarrow}(Y_t)\right| =  \left|\lambda_j^{\downarrow}(X)\right| - \left|\lambda_j^{\downarrow}(Y_t)\right|
\peso{for every} 1\le j \le n \ .
$$
Making the same computation by replacing $X$ by $Y_s$ for any $s\in (t\coma 1]$ 
(by the minimality of $\alpha$ and the fact that $\|Y_s\|_{\infty}<\pi$, it follows that $d(I\coma \alpha(s)\,) = \|Y_s\|_1$), we deduce that the maps $t\mapsto \left|\lambda_j^{\downarrow}(Y_t)\right|$ are increasing 
and that $\lambda_j(X), \lambda_j(Y_t) \ge 0$ (they have the same sign) 
for every $1\le j \le n$ and every $t \in [0, 1]$.
%
Therefore 
the maps $t\mapsto \lambda_j^{\downarrow}(Y_t)$ are monotone and minimal in $\R$ and also 
the curves  $\la_j(\alpha (t)) = e^{i\, \la_j\da(Y_t\,)}$ 
are minimal in $\C$ for trace norm (the second part
follows using Theorem \ref{corta Un} for $n=1$).
\edem

\newpage
\section{Intermediate curves and geometry of midpoints}

In this section we will prove that the set of midpoints is geodesically convex in all the previous contexts. Actually, we will use the same idea for all this cases. First we prove that the function
$$
s \mapsto d(I, \gamma(s)),
$$
where $\gamma(s)$ is the geodesic joining two matrices (in some of this spaces for some metric), is convex. Then we will use this fact to prove that the set of midpoints is geodesically convex.

\subsection*{The unitary group}

Let $U,\, V\in\matu$ and $t\in(0,1)$. In this section we will study the sets of intermediate points for the spectral norm:
$$
\eme_t(U,V)=\{W\in\matu: d_\infty(U,W)=t\,d_\infty(U,V)\  \mbox{and}\  d_\infty(W,V)=(1-t)\,d_\infty(U,V) \}.
$$
In the case of positive matrices, the corresponding sets of intermediate points were proved to be geodesically convex (see \cite{L1}). We can not expect a similar result with full generality in our setting. Indeed, for instance, in the one dimensional case, the set of intermediate points may be even disconnected. In the case of dimension greater than one, there is also a particular case, where the geodesic convexity is not true:

\begin{exa}
	Take $U=I$ and $V=-I$, and consider the set of intermediate points $\eme_t(I,-I)$ for some $t\in(0,1/2)$ (for $t\in(1/2,1)$ is similar). Then
	$$
	W^+=\begin{pmatrix}
	e^{it\pi}&0\\
	0& e^{it\pi}
	\end{pmatrix}\peso{and}
	W^-=\begin{pmatrix}
	e^{-it\pi}&0\\
	0& e^{-it\pi}
	\end{pmatrix}.
	$$
	belongs to $\eme_t(I,-I)$. However, it is not difficult to see that $I$ is the midpoint of the geodesic (and unique minimal curve) that joints $W^+$ with $W^-$. \EOE
\end{exa}

\medskip

First of all, we will give a proof of the following characterization of the sets $\eme_t(U,V)$.

\begin{lem}
	Given $U,W\in\matu$ and $t\in(0,1)$, then
	\begin{align*}
	M_t(U, V) & = \left\{V \in \mathcal{U}_n :  \ \mbox{exists}  \ \gamma : [0, 1] \rightarrow \cU(n) \ \mbox{a minimal curve} \ : \ \gamma(t) = V \right\}.
	\end{align*}
\end{lem}
\bdem
If there exist a minimal curve $\gamma$ such that $\gamma(t) = W$ then, by Corollary \ref{lineal}
\begin{align*}
d_\infty(U, W) & = d_\infty(\gamma(0), \gamma(t)) = td_\infty(U,V) \\
d_\infty(W, V) & = d_\infty(\gamma(t), \gamma(1)) = (1- t)d_\infty(U,V).
\end{align*}
This proves one inclusion. Conversely, if $W\in M_t(U, W)$, take a minimal curve $\beta_1$ joining $U$ with $W$ and $\beta_2$ a minimal curve joining $W$ with $V$. Then
\begin{align*}
\mbox{L}(\gamma) &= \mbox{L}(\beta_1) + \mbox{L}(\beta_2) \\
&= d_\infty(U, W) + d_\infty(W, V) \\
&=  td_\infty(U,W) + (1 - t)d_\infty(W, V) =d_\infty(U, V).
\end{align*}
If $\gamma = \beta_1 * \beta_2$ is the concatenation of these two curves parametrized with constant speed in $[0,1]$, then $W=\gamma(t)$.
\edem

\bigskip

We are going to prove that if $W_0, W_1 \in M_{t}(U, V)$ then $\beta(s) \in M_{t}(U, V)$, for all $s \in [0, 1]$, where $\beta(s)$ is a geodesic with the condition $d_\infty(U, V) < \pi/2$. For that propose we will use the following theorem.

\begin{teo}[See \cite{AL}, Thm 2.8]
	
	If $d_\infty(W, \beta(s)) < \pi/2$ for all $s \in [0, 1]$, where $\beta(s)$ is a geodesic, then 
	$$
	s \mapsto d_\infty(W, \beta(s))
	$$
	is convex.
\end{teo}

If $d_\infty(U, V) < \pi/2$, given $W_0, W_1 \in M_{t}(U, V)$, 
\begin{align*}
d_\infty(U, \beta(s)) & \leq d_\infty(U, W_j) + \frac{d_\infty(W_0, W_1)}{2} \\
& \leq d_\infty(U, W_j) + \frac{d_\infty(W_0, V)}{2} +   \frac{d_\infty(W_1, V)}{2} \\ 
& = td_\infty(U, V) + \frac{(1-t)d_\infty(U, V)}{2} +   \frac{(1-t)d_\infty(U, V)}{2} = d_\infty(U, V) < \frac{\pi}{2}
\end{align*}
taking a $j =0, 1$ convenient.

\begin{figure}[h!]
	\centering
	\includegraphics[width=5cm]{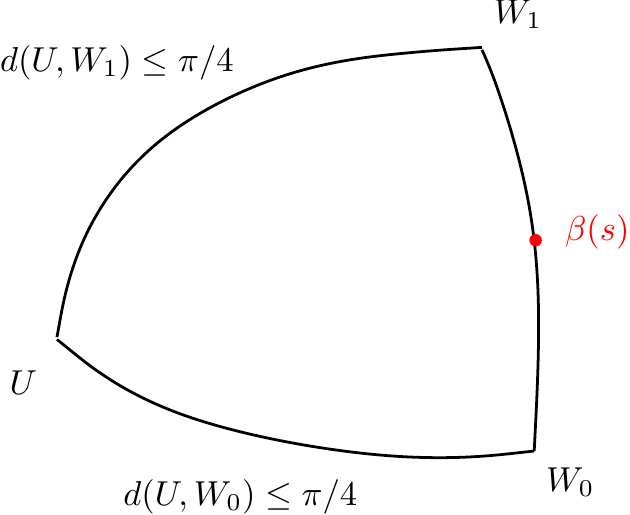} 
\end{figure} 

\begin{pro}\label{convexity of midpoints}
	Given $U,V\in\matu$ such that $d_\infty(U,V)<\pi/2$, and $t\in(0,1)$, the set
	$$
	\eme_t(U,V)=\{W\in\matu: d_\infty(U,W)=t\,d_\infty(U,V)\  \mbox{and}\  d_\infty(W,V)=(1-t)\,d_\infty(U,V) \}
	$$
	is geodesically convex.
\end{pro}

\bdem
Note that, if $\beta(s)$ is the geodesic joining $W_o, W_1 \in \eme_t(U,V)$, i.e. $W_0 = \beta(0), W_1 = \beta(1)$, then 
\begin{align*}
d_\infty(U, \beta(s)) & \leq s d_\infty(U, W_1) + (1 - s)d_\infty(U, W_0) \\
& = s t d_\infty(U, V) + (1 - s)t d_\infty(U, V) = t d_\infty(U, V); \\
d_\infty(V, \beta(s)) & \leq s d_\infty(V, W_1) + (1 - s)d_\infty(V, W_0) \\
& = s (1-t) d_\infty(U, V) + (1 - s)(1-t) d_\infty(U, V) = (1-t) d_\infty(U, V)
\end{align*}
Therefore
$$
d_\infty(U, V) \leq d_\infty(U, \beta(s)) + d_\infty(V, \beta(s)) \leq d_\infty(U, V).
$$
That is,
\begin{align*}
d_\infty(U, \beta(s)) &= td_\infty(U, V) \\
d_\infty(V, \beta(s)) &= (1-t)d_\infty(U, V).
\end{align*}
Then $\beta(s) \in M_{t}(U, V)$.
\edem

\subsection*{The Grassmannian}

As before, our previous results about the unitary group can by applied to the Grassmann manifolds. 

\begin{pro}\label{convexity of midpoints2}
	Given $P, Q \in \mathcal{G}_n$ such that $\left\|P - Q\right\|_{\infty} < 1/\sqrt{2}$, and $t\in(0,1)$, the set
	$$
	\eme_t(P,Q)=\{R\in\matu: d_\infty(P,R)=t\,d_\infty(P,Q)\  \mbox{and}\  d_\infty(R,Q)=(1-t)\,d_\infty(P,Q) \}
	$$
	is geodesically convex.
\end{pro}

\bdem
By Proposition \ref{convexity of midpoints}, if $d(S_P, S_Q) < \pi/2$ then the set $M_t(S_P, S_Q)$ is geodesically convex. As $d(S_P, S_Q) < \pi/2$ is equivalently to $\left\|S_P - S_Q\right\|_{\infty} < \sqrt{2}$ then
$$
\left\|P - Q\right\|_{\infty}  = \frac{1}{2}\left\|2P - 2Q\right\|_{\infty}  = \frac{1}{2}\left\|(2P - 1) - (2Q - 1)\right\|_{\infty}  = \frac{1}{2}\left\|S_P - S_Q\right\|_{\infty}  < \frac{1}{2} \sqrt{2}.
$$
implies that $M_t(P, Q)$ is geodesically convex.
\edem

\subsection*{The spaces 
of positive invertible matrices and 
of Hermitian matrices}

Let $P\in\matppos$ and $t\in(0,1)$. In what follows we fix $N$ a UIN on $\mat$ and we 
 study the sets of intermediate points 
$$
\eme_t(I,P)=\{W\in\matppos :  d_N (I,W)=t\, d_N (I,P)\  \mbox{and}\   d_N (W,P)=(1-t)\, d_N (I,P) \}.
$$
If $N$ is the trace norm, the corresponding sets of intermediate points were proved to be geodesically convex (see \cite{L1}). 
Here we present a proof which is much simpler and is valid for every UIN. It is based in  the same idea as the one given for 
$\matu$, 
and the fact that the map 
$$
s \mapsto  d_N (I, \gamma(s)),
$$
where $\gamma(s)$ is the geodesic joining two positive matrices, is convex. This fact is 
known, but we present a short proof for completeness. 
\begin{teo}\label{58}
Fix $N$ a UIN in $\mat$. 
If $\gamma(s) = A^{1/2}(A^{-1/2}BA^{-1/2})^sA^{1/2}$ denote the geodesic joining $A, B \in \matppos$, then the map
	$$
	[0\coma 1] \ni s \longmapsto  d_N (I, \gamma(s))
	$$
	is convex.
\end{teo}

\bdem
By Araki inequality (see \cite{Araki}), it is easy to see that 
$$
d_N (C^t \coma D^t)  = N\big(\log \,( C^{t/2} D^{-t} C^{t/2} ) \,\big)\le 
t  \, N\big( \log \,( C^{1/2} D^{-1} C^{1/2} ) \,\big)= 
t  \, d_N (C \coma D) \ ,
$$ 
for every $t\in [0\coma 1] $ and every pair $C\coma D \in \matppos$. Then, for $s\in  [0\coma 1] $, 
\begin{align*}
 d_N (I, \gamma(s)) & \leq  d_N (I, A^{1-s}) +  d_N (A^{1-s}, \gamma(s)) \\ 
 & =  d_N (I, A^{1-s}) +  d_N (A^{-s}, (A^{-1/2}BA^{-1/2})^s) \\ & \leq 
 (1-s) \, d_N (I, A) + s \, d_N (A^{-1}, A^{-1/2}BA^{-1/2}) \\ & = (1-s) \, d_N (I, A) + s \, d_N (I, B).
\end{align*}
It is easy to see that this fact implies the convexity of our map. 
\edem

\begin{pro}\label{57}
Let $P\in\matppos$ and $t\in(0,1)$. Then the set
$$
\eme_t(I,P)=\{W\in\matppos :  d_N (I,W)=t\, d_N (I,P)\  \mbox{and}\   d_N (W,P)=(1-t)\, d_N (I,P) \}
$$
is geodesically convex.
\end{pro}

\bdem
Let $\gamma(s)= W_0^{1/2}(W_0^{-1/2}W_1W_0^{-1/2})^sW_0^{1/2}$ be the geodesic joining $W_0$ and $ W_1 \in \eme_t(U,V)$. Then 
\begin{align*}
 d_N (I, \gamma(s)) & \leq (1-s) d_N (I, W_0) + s d_N (I, W_1) \\
& = (1-s) t  d_N (I, P) + st  d_N (I, P) = t  d_N (I, P); \\
 d_N (P, \gamma(s)) & \leq (1-s)  d_N (P, W_0) + s d_N (P, W_1) \\
& = (1-s)(1-t)  d_N (I, P) + s(1-t)  d_N (I, P) = (1-t)  d_N (I, P).
\end{align*}
Hence $
 d_N (I, P) \leq  d_N (I, \gamma(s)) +  d_N (P, \gamma(s)) \leq  d_N (I, P)
$, so that 
$$
 d_N (I, \gamma(s)) = t d_N (I, P) \peso{and} 
 d_N (P, \gamma(s)) = (1-t) d_N (I, P).	
$$
Then $\gamma(s) \in M_{t}(I, P)$ for all $s \in [0, 1]$.
\edem

%
%
%
%
%
%

\begin{rem}
If the norm $N$ is strictly convex then the geodesic joining two matrices in $\matppos$ or $\matsa$ is unique, 
and the set $\eme_t(I,P)$ is a single point. 
Nevertheless, Proposition \ref{57} is interesting 
for several important UIN's, such as the trace, spectral and other Ky Fan norms. 
	\EOE
\end{rem}

\section*{Acknowledgements:}
This work was supported by Consejo Nacional de Investigaciones Cient\'\i ficas y T\'ecnicas-Argentina (PIP-152), Agencia Nacional de Promoci\'on de Ciencia y Tecnolog\'\i a-Argentina (PICT 2015-1505), Universidad Nacional de La Plata-Argentina (UNLP-11X585) and  Ministerio de Econom\'ia y Competitividad-Espa\~na (MTM2016-75196-P).

\end{document}